\documentclass[12pt]{article}

\usepackage{amssymb}
\usepackage{amsmath}
\usepackage{cite}
\usepackage[arrow, matrix, curve]{xy}
\usepackage{bbm}

\usepackage[dvips]{graphicx}
\usepackage{xcolor}

\begin{document}

\renewcommand{\citeleft}{{\rm [}}
\renewcommand{\citeright}{{\rm ]}}
\renewcommand{\citepunct}{{\rm,\ }}
\renewcommand{\citemid}{{\rm,\ }}

\newcounter{abschnitt}
\newtheorem{satz}{Theorem}
\newtheorem{coro}[satz]{Corollary}
\newtheorem{theorem}{Theorem}[abschnitt]
\newtheorem{koro}[theorem]{Corollary}
\newtheorem{prop}[theorem]{Proposition}
\newtheorem{lem}[theorem]{Lemma}
\newtheorem{expls}[theorem]{Examples}

\newcommand{\mres}{\mathbin{\vrule height 1.6ex depth 0pt width 0.11ex \vrule height 0.11ex depth 0pt width 1ex}}

\renewenvironment{quote}{\list{}{\leftmargin=0.62in\rightmargin=0.62in}\item[]}{\endlist}

\newcounter{saveeqn}
\newcommand{\alpheqn}{\setcounter{saveeqn}{\value{abschnitt}}
\renewcommand{\theequation}{\mbox{\arabic{saveeqn}.\arabic{equation}}}}
\newcommand{\reseteqn}{\setcounter{equation}{0}
\renewcommand{\theequation}{\arabic{equation}}}

\hyphenpenalty=9000

\sloppy

\phantom{a}

\vspace{-1.5cm}

\begin{center}
\begin{Large} {\bf Blaschke--Santal\'o Inequalities for\\ Minkowski and Asplund Endomorphisms} \\[0.7cm] \end{Large}

\begin{large} Georg C. Hofst\"atter and Franz E.\ Schuster\end{large}
\end{center}

\vspace{-0.9cm}

\begin{quote}
\footnotesize{ \vskip 1cm \noindent {\bf Abstract.}
It is shown that each monotone Minkowski endomorphism of convex bodies gives rise
to an isoperimetric inequality which directly implies the classical Urysohn inequality. Among this large family of new inequalities, the only affine invariant one -- the Blaschke--Santal\'o inequality -- turns out to be the strongest one. A further extension of these inequalities to merely weakly monotone Minkowski endomorphisms is proven to be impossible. Moreover, for functional analogues of monotone Minkowski endomorphisms, a family of analytic inequalities for log-concave functions is established which generalizes the functional Blaschke--Santal\'o inequality.

}
\end{quote}

\vspace{0.6cm}

\centerline{\large{\bf{ \setcounter{abschnitt}{1}
\arabic{abschnitt}. Introduction}}}

\alpheqn

\vspace{0.5cm}

One of the most widely known and fundamental affine isoperimetric inequalities is the \emph{Blaschke--Santal\'o inequality}, roughly stating that the volume product of polar reciprocal convex bodies is maximized by ellipsoids. More precisely, let $K \subseteq \mathbb{R}^n$ be a convex body (that is, a compact, convex set) with non-empty interior and recall that
$K^z = \{x \in \mathbb{R}^n: x \cdot y \leq 1 \mbox{ for all } y \in K - z\}$ is the \emph{polar body} of $K$ with respect to $z \in \mathrm{int}\,K$. Denoting by $|A|$ the volume of a Borel set $A \subseteq \mathbb{R}^n$, the \emph{Santal\'o point} can be defined as the unique point $\mathbf{s} = \mathbf{s}(K) \in \mathrm{int}\,K$, for which
$|K^\mathbf{s}| = \min \{|K^z|: z \in \mathrm{int}\,K  \}$. The Blaschke--Santal\'o inequality then states that
\begin{equation} \label{BSinequ}
|K||K^\mathbf{s}| \leq |B^n|^2
\end{equation}
with equality if and only if $K$ is an ellipsoid (that is, an affine image of the Euclidean unit ball $B^n$).
Initial proofs of (\ref{BSinequ}) were given in the first half of the previous century by Blaschke for $n \leq 3$ and Santal\'o for all $n \geq 2$, while the equality conditions
were completely settled only in 1985 by Petty. In subsequent years, simplified proofs, including the equality cases, were obtained (see, e.g., \textbf{\cite{FradMeyer2007, MeyerPajor1989, MeyerPajor1990, MeyerReisner2006}}) and it remained an active focus of research due to the evolving understanding of its impact (see
\textbf{\cite{alexfradzvav2019, boroczky2010, corderetal2015, ivaki2015a, LutwakZhang1997, MeyerWerner1998, Schneider1982}} and the references therein).

Affine invariant inequalities are often more powerful than related inequalities that are merely invariant under Euclidean rigid motions. This becomes particularly striking for the Blaschke--Santal\'o inequality which considerably strengthens and directly implies the classical \emph{Urysohn inequality} (as first observed by Lutwak \textbf{\cite{Lutwak1975}}). The latter is the following basic relation between the mean width $w(K)$ of a convex body $K \subseteq \mathbb{R}^n$ (see Section 2 for definition) with non-empty interior and its volume,
\begin{equation} \label{Uryinequ}
|K| \leq \left (\frac{w(K)}{2}\right )^n |B^n|
\end{equation}
with equality if and only if $K$ is a ball.

\pagebreak

Another affine isoperimetric inequality that plays a special role in this paper coincides for origin-symmetric bodies with the Blaschke--Santal\'o inequality but is in general weaker than (\ref{BSinequ}). In order to state it, let $\Delta K = \frac{1}{2}(-K + K)$ denote the \emph{central symmetral} of a convex body $K \subseteq \mathbb{R}^n$. If $K$ has non-empty interior, then (\ref{BSinequ}), combined with the Brunn--Minkowski inequality, implies that
\begin{equation} \label{diffbodinequ}
|K||\Delta^{\circ}K| \leq |B^n|^2
\end{equation}
with equality if and only if $K$ is an ellipsoid. Here, $\Delta^\circ K$ is the polar body of $\Delta K$ with respect to the origin. The central symmetrization $\Delta$ has long been a useful tool \linebreak in the Brunn--Minkowski theory (see, e.g., \textbf{\cite[\textnormal{Chapter 3.2}]{gardner2ed}} and \textbf{\cite[\textnormal{Chapter~10.1}]{schneider93}}). As a continuous operator on the space $\mathcal{K}^n$ of convex bodies in $\mathbb{R}^n$ endowed with the Hausdorff metric, the importance of $\Delta$ stems from its \emph{Minkowski additivity}  (that is, $\Delta(K + L) = \Delta K + \Delta L$ for all $K, L \in \mathcal{K}^n$) and compatibility with affine transformations. These are characterizing properties, as the following result shows.

\vspace{0.2cm}

\noindent {\bf Theorem} {(Schneider \textbf{\cite{Schneider1974c}})}
\emph{A continuous map $\Phi: \mathcal{K}^n \rightarrow \mathcal{K}^n$ is a translation invariant Minkowski additive map such that $\Phi(AK) = A\Phi K$ for every $K \in \mathcal{K}^n$ and $A \in \mathrm{GL}(n)$ if and only if $\Phi = c\,\Delta$ for some $c \geq 0$.}

\vspace{0.2cm}

This theorem was a byproduct of a more general, systematic study of Minkowski additive operators on $\mathcal{K}^n$, initiated about 50 years ago by Schneider \textbf{\cite{Schneider1974, Schneider1974b, Schneider1974c}}.
Since then, and up to now, the main focus thereby has been on maps that also commute with $\mathrm{SO}(n)$ transforms (see \textbf{\cite{Abardiaetal2018, Dorrek2017b, Kiderlen2006, Schu06a, Schu09, SchuWan16}}). As such maps are automatically compatible with translations (see, e.g., \textbf{\cite[\textnormal{Section 2.3}]{Kiderlen2006}}), they are often assumed w.l.o.g.\ to be translation invariant, leading to the following central definition.

\vspace{0.2cm}

\noindent {\bf Definition} \emph{A continuous map $\Phi: \mathcal{K}^n \rightarrow \mathcal{K}^n$ is a \emph{Minkowski endomorphism} if
$\Phi$ is Minkowski additive, translation invariant, and commutes with $\mathrm{SO}(n)$ transforms. The \emph{trivial} Minkowski endomorphism maps every convex body to the origin.}

\vspace{0.2cm}

Much of this article is motivated by the observation that (\ref{Uryinequ}) and (\ref{diffbodinequ}) can be cast as  volume estimates for polar Minkowski endomorphisms. Another prominent such example was established by Lutwak \textbf{\cite{Lutwak1993}} for polar projection bodies of order~1,
\begin{equation} \label{Pi1inequ}
|K||\Pi_1^\circ K| \leq |B^n|^2
\end{equation}
with equality if and only if $K$ is a ball. Recalling that each $K \in \mathcal{K}^n$ is determined by its support function $h(K,u) = \max\{u \cdot x: x \in K\}$ for $u \in \mathbb{S}^{n-1}$, $\Pi_1 K$ can be defined by $h(\Pi_1K,u) = c_n w(K|u^\perp)$, where $K|u^\perp$ denotes the orthogonal projection of $K$ onto $u^\perp$ and $c_n \in \mathbb{R}$ is chosen such that $\Pi_1B^n = B^n$.

The natural question to what degree inequalities (\ref{Uryinequ}), (\ref{diffbodinequ}), and (\ref{Pi1inequ}) can be unified, was first asked by Lutwak. A partial answer was given in \textbf{\cite{bergschu2020}}, deduced from results in \textbf{\cite{HaberlSchu2019}}, where (\ref{Uryinequ}) and (\ref{Pi1inequ}) were identified as part of a larger family of inequalities for a subcone of Minkowski endomorphisms which are \emph{monotone}, that is, $K \subseteq L$ implies $\Phi K \subseteq \Phi L$ for all $K, L \in \mathcal{K}^n$. For a more precise statement we require the following classification of monotone Minkowski endomorphisms.

\pagebreak

% \vspace{0.3cm}

\noindent {\bf Theorem} {(Kiderlen \textbf{\cite{Kiderlen2006}})}
\emph{A map $\Phi: \mathcal{K}^n \rightarrow \mathcal{K}^n$ is a monotone Minkowski endo\-morphism if and only if there exists a non-negative $\mathrm{SO}(n-1)$ invariant measure $\mu$ on $\mathbb{S}^{n-1}$ with center of mass at the origin such that
\begin{equation} \label{kidmainformular}
h(\Phi K,\cdot) = h(K,\cdot) \ast \mu
\end{equation}
for every $K \in \mathcal{K}^n$. Moreover, the measure $\mu$ is uniquely determined by $\Phi$.}

\vspace{0.3cm}

The convolution of functions and measures on $\mathbb{S}^{n-1}$ used in (\ref{kidmainformular}) is induced from $\mathrm{SO}(n)$ by identifying $\mathbb{S}^{n-1}$ with the homogeneous space $\mathrm{SO}(n)/\mathrm{SO}(n-1)$ (see Section 2 for details). Note that we assume all measures to be finite Borel measures.

In \textbf{\cite{bergschu2020}}, (\ref{Uryinequ}) and (\ref{Pi1inequ}) were generalized to monotone Minkowski endomorphisms generated by area measures of order one of zonoids (see Section 2). As a first main result, we generalize these inequalities from \textbf{\cite{bergschu2020}} to \emph{all} monotone Minkowski endomorphisms $\Phi$. Throughout, we always assume that $n \geq 3$.

\begin{satz} \label{thm:ineqVolProdMonotone}
Suppose that $\Phi: \mathcal{K}^n \rightarrow \mathcal{K}^n$ is a monotone non-trivial Minkowski endomorphism. Among $K \in \mathcal{K}^n$ with non-empty interior the volume product
\[|K||\Phi^\circ K|  \]
is maximized by Euclidean balls. If $\Phi = c\,\Delta$ for some $c > 0$, then $K$ is a maximizer if and only if it is an ellipsoid. Otherwise, Euclidean balls are the only maximizers.
\end{satz}

Let us emphasize that Theorem \ref{thm:ineqVolProdMonotone} not only includes inequalities (\ref{Uryinequ}), (\ref{diffbodinequ}), and (\ref{Pi1inequ}) as special cases, but provides an extension of the isoperimetric inequalities from \textbf{\cite{bergschu2020}} from a nowhere dense set of Minkowski endomorphisms to all monotone ones. Whereas the proof of Theorem \ref{thm:ineqVolProdMonotone} does not require any results from \textbf{\cite{HaberlSchu2019}}, our approach is very much inspired by techniques from \textbf{\cite{HaberlSchu2019}} and relies on Kiderlen's classification of monotone Minkowski endomorphisms.

While it was long known that not all Minkowski endomorphisms are monotone, a conjecture that they are all \emph{weakly monotone} (see Section 3 for details) was disproved by Dorrek \textbf{\cite{Dorrek2017b}} only recently.
We will show in Section 4 that Theorem \ref{thm:ineqVolProdMonotone} is essentially the best possible, in the sense that a further extension to all merely weakly monotone endomorphisms is in general impossible.

By Schneider's above characterization of the  map $\Delta$, inequality (\ref{diffbodinequ}) is the only affine invariant one among the family of isoperimetric inequalities provided by Theorem \ref{thm:ineqVolProdMonotone}. With our second main result, we show that all these inequalities can be deduced from the Blaschke--Santal\'o inequality. In particular, among inequalities for even Minkowski endomorphisms, (\ref{diffbodinequ}) is the strongest member of the inequalities from Theorem \ref{thm:ineqVolProdMonotone}. This is in contrast to the volume estimates obtained in \textbf{\cite{bergschu2020}}, among which (\ref{Pi1inequ}) was the strongest one, since (\ref{diffbodinequ}) was not included.  Finally, we prove that each of the inequalities of Theorem \ref{thm:ineqVolProdMonotone} is stronger and directly implies the Urysohn inequality (\ref{Uryinequ}).

\pagebreak

\begin{satz} \label{thm:ineqUrysohnMonotone}
If $\Phi: \mathcal{K}^n \rightarrow \mathcal{K}^n$ is a monotone Minkowski endomorphism such that $\Phi B^n = B^n$ and $K \in \mathcal{K}^n$ has non-empty interior, then
\begin{align} \label{thm2inequ}
|B^n|\left(\frac{w(K)}{2} \right)^{-n} \leq |\Phi^\circ K| \leq | K^{\mathbf{s}}|.
\end{align}
There is equality in the left hand inequality if and only if $\Phi K$ is a Euclidean ball. Equality in the right hand inequality holds if and only if
$K$ is centrally-symmetric and $\Phi = \Delta$ or if $K$ is a Euclidean ball.
\end{satz}

Note that the right inequality of Theorem \ref{thm:ineqUrysohnMonotone} (and its equality case) combined with the Blaschke--Santal\'o inequality implies Theorem \ref{thm:ineqVolProdMonotone}. In Section 4, we will therefore first prove Theorem \ref{thm:ineqUrysohnMonotone} and then deduce Theorem \ref{thm:ineqVolProdMonotone} as a consequence.

\vspace{0.15cm}

A second focus of this article concerns the continuing effort to extend notions and results from convex geometry to the class of log-concave functions, that is, all $f: \mathbb{R}^n \rightarrow [0,\infty)$ of the form $f = e^{-\varphi}$ for some convex $\varphi: \mathbb{R}^n \rightarrow (-\infty,\infty]$. The most basic such notions are Minkowski addition and scalar multiplication which can be naturally extended as follows. For log-concave $f$ and $g$ and $\lambda > 0$, let
\[(f \star g)(x) = \sup_{x_1 + x_2 = x} f(x_1)g(x_2), \qquad (\lambda \cdot f)(x) = f\left (\mbox{$\frac{x}{\lambda}$}\right)^{\lambda}.  \]
Then $f \star g$ is called the \emph{Asplund sum} (or sup-convolution) of $f$ and $g$ (see, e.g., \textbf{\cite{gardkid2018}}). \linebreak
While the above definitions imply that $\mathbbm{1}_{K} \star \mathbbm{1}_{L} = \mathbbm{1}_{K+L}$ and $\lambda \cdot \mathbbm{1}_K = \mathbbm{1}_{\lambda K}$ for all indicators of $K, L \in \mathcal{K}^n$ and $\lambda > 0$, in general it is possible that $f \star g$ attains the value $+\infty$ and, thus, is no longer log-concave. Moreover, the standard regularity assumption of upper semi-continuity of log-concave functions need not be preserved under Asplund addition (cf.\ \textbf{\cite[\textnormal{p.\ 517}]{schneider93}}). One frequently used possibility to overcome these issues is to work with the space $\mathrm{LC}_{\mathrm{c}}(\mathbb{R}^n)$ of all proper log-concave functions which are upper semi-continuous and coercive. Here, $f$ is called \emph{proper} if it is not identically $0$ and it is \emph{coercive} if $\lim_{\|x\|\rightarrow \infty} f(x) = 0$.
We furthermore endow $\mathrm{LC}_{\mathrm{c}}(\mathbb{R}^n)$ with the topology induced by epi-convergence (see Section 2).

As a seminal inequality for log-concave functions, we first mention the celebrated Pr\'ekopa--Leindler inequality which is universally recognized as the functional form of the Brunn--Minkowski inequality (see, e.g., \textbf{\cite[\textnormal{Section 7.1}]{schneider93}}). A functional version of the Blaschke--Santal\'o inequality, discovered later by Ball \textbf{\cite{Ball1986}}, is of special importance for us. Before stating it, recall that for log-concave $f: \mathbb{R}^n \rightarrow [0,\infty)$, the \emph{polar function} of $f$ can be defined (following \textbf{\cite{Artstein2004}}) by
\begin{equation} \label{defpolarf}
f^{\circ} = e^{-\mathcal{L}(-\log f)},
\end{equation}
where $\mathcal{L}$ denotes the classical Legendre transform (see Section 2). If $f$ is in addition \emph{even} and $0 < \int_{\mathbb{R}^n}\! f dx< \infty$, then Ball's \emph{functional Blaschke--Santal\'o inequality} reads
\begin{equation} \label{ballfunctBSinequ}
\int_{\mathbb{R}^n} f(x)\,dx  \int_{\mathbb{R}^n} f^\circ(x)\,dx \leq (2\pi)^n
\end{equation}
with equality if and only if $f$ is a Gaussian.

\pagebreak

The cases for equality in (\ref{ballfunctBSinequ}) were settled by  Artstein-Avidan, Klartag and Milman~\textbf{\cite{Artstein2004}}, who also established a far-reaching extension of (\ref{ballfunctBSinequ}) to not necessarily even functions (cf.\ Section 2), that has sparked a great deal of research interest in recent years, see \textbf{\cite{Artstein2015b, BartheBoerFrad2014, ColFrag2013, Cordero2015b, FradMeyer2007, Haddadetal2020, KolWerner2021, Lehec2009, Lehec2009a, Rotem2014, Rotem2020}}.

As noted by Rotem \textbf{\cite{Rotem2012}}, the functional Blaschke--Santal\'o inequality implies an analogue of Urysohn's inequality for log-concave functions, a result first obtained by Klartag and Milman \textbf{\cite{KlartagMilman2005}} by other means. It can be conveniently formulated with the help of the support function of a log-concave $f : \mathbb{R}^n \rightarrow [0,\infty)$ which is, following \textbf{\cite{Artstein2010}}, defined by $h(f,\cdot) = \mathcal{L}(-\log f)$. If additionally $\int_{\mathbb{R}^n} f\,dx = (2\pi)^{n/2}$ and $\gamma_n$ is the standard Gaussian measure on $\mathbb{R}^n$, the functional analogue of Urysohn's inequality states that
\begin{equation} \label{functurysohn}
\frac{2}{n} \int_{\mathbb{R}^n} h(f,x)\,d\gamma_n(x) \geq 1
\end{equation}
with equality if and only if $(2\pi)^{-n/2}f$ is a translation of the standard Gaussian. Note that if $f = \mathbbm{1}_K$ for some $K \in \mathcal{K}^n$, the left hand side is proportional to $w(K)$. However, the \emph{sharp} geometric Urysohn inequality cannot be recovered from (\ref{functurysohn}).

Recalling that $\mathrm{SO}(n)$ acts naturally on $\mathrm{LC_c}(\mathbb{R}^n)$, specifically, $(\vartheta f)(x) = f(\vartheta^{-1}x)$ for $\vartheta \in \mathrm{SO}(n)$ and $f \in \mathrm{LC_c}(\mathbb{R}^n)$, we can now introduce `functional Minkowski endomorphisms' on $\mathrm{LC_c}(\mathbb{R}^n)$ as follows.

\vspace{0.3cm}

\noindent {\bf Definition} \emph{A continuous map $\Psi: \mathrm{LC_c}(\mathbb{R}^n) \rightarrow \mathrm{LC_c}(\mathbb{R}^n)$ is an \emph{Asplund endomor- phism} if
$\Psi$ is Asplund additive, translation invariant, and commutes with the $\mathrm{SO}(n)$ action. The \emph{trivial} Asplund endomorphism maps every function to the indicator of the origin. (Note that we do not identify functions that coincide almost everywhere.)}

\vspace{0.3cm}

In Section 3, we discuss the problem of establishing an analogue of Kiderlen's characterization of monotone Minkowski endomorphisms. Here, we use the latter as motivation to define the following rich class of Asplund endomorphisms which are \emph{monotone}, that is, $f \leq g$ implies $\Psi f \leq \Psi g$ for all $f,g \in \mathrm{LC_c}(\mathbb{R}^n)$.

\begin{satz} \label{def:phimuFirst}
Each non-negative $\mathrm{SO}(n-1)$ invariant measure $\mu$ on $\mathbb{S}^{n-1}$ with center of mass at the origin induces a monotone Asplund endomorphism $\Psi_{\mu}$ by
\begin{equation} \label{thm3introform}
h(\Psi_{\mu} f,\cdot) = h(f,\cdot) \circledast \mu 
\end{equation}
for $f  \in \mathrm{LC_c}(\mathbb{R}^n)$. Moreover, the measure $\mu$ is uniquely determined by $\Psi_{\mu}$.
\end{satz}

For the definition of the convolution $\circledast$ of the convex function $h(f,\cdot)$ on $\mathbb{R}^n$ with the measure $\mu$, we refer to Section 2, where it will also become evident that (\ref{thm3introform}) defines a mapping from $\mathrm{LC_c}(\mathbb{R}^n)$ to itself. Let us emphasize that $\Psi_{\mu} \mathbbm{1}_K = \mathbbm{1}_{\Phi_{\mu} K}$ for every $K \in \mathcal{K}^n$, where $\Phi_{\mu}$ is the monotone Minkowski endomorphism defined by (\ref{kidmainformular}). In this sense, the Asplund endomorphisms $\Psi_{\mu}$ extend the class of all monotone Minkowski endomorphisms to $\mathrm{LC_c}(\mathbb{R}^n)$. As our next main result, we prove a functional analogue of Theorem \ref{thm:ineqVolProdMonotone} for the monotone Asplund endomorphisms defined by Theorem \ref{def:phimuFirst}.

\pagebreak

\begin{satz} \label{thm:ineqVolProdFunc}
Let $\mu$ be an $\mathrm{SO}(n-1)$ invariant probability measure on $\mathbb{S}^{n-1}$ with center of mass at the origin. If $f\in \mathrm{LC_c}(\mathbb{R}^n)$ such that $\int_{\mathbb{R}^n} f\,dx > 0$, then
 \begin{align} \label{functinequAsplendos}
  \int_{\mathbb{R}^n} f(x)\,dx  \int_{\mathbb{R}^n}(\Psi_{\mu} f)^\circ(x)\, dx \leq (2\pi)^n.
 \end{align}
If $\mu$ is discrete, there is equality if and only if $f$ is a Gaussian. Otherwise, equality holds if and only if $f$ is proportional to a translation of the standard Gaussian.
\end{satz}

Note that the additional normalization of $\mu$ in Theorem \ref{thm:ineqVolProdFunc} is critical due to the non-homogeneity of the integral as well as the functional polarity with respect to Asplund scalar multiplication.

\vspace{0.3cm}

\noindent {\bf Examples:}

\begin{enumerate}
\item[(a)] In 2006, Colesanti \textbf{\cite{colesanti2006}} introduced the (Asplund) \emph{difference function} $\Delta_{\star}f$ of a log-concave function $f \in \mathrm{LC_c}(\mathbb{R}^n)$ by
$\Delta_{\star} f = \frac{1}{2} \cdot f \star \frac{1}{2} \cdot \overline{f}$, where $\overline{f}(x) = f(-x)$. \linebreak By taking $\mu$ to be the even discrete probability measure $\frac{1}{2}\delta_{\bar{e}} + \frac{1}{2}\delta_{-\bar{e}}$, concentrated on the stabilizer $\bar{e}\in \mathbb{S}^{n-1}$ of $\mathrm{SO}(n-1)$ and its antipodal, Theorem \ref{thm:ineqVolProdFunc} reduces to a functional analogue of (\ref{diffbodinequ}),
\begin{equation} \label{deltaasplinequ}
 \int_{\mathbb{R}^n} f(x)\,dx  \int_{\mathbb{R}^n}(\Delta_{\star} f)^\circ(x)\, dx \leq (2\pi)^n
\end{equation}
which can also be deduced from Ball's functional Blaschke--Santal\'o inequality (\ref{ballfunctBSinequ}) and the Pr\'ekopa--Leindler inequality.
Clearly, for even $f$, (\ref{deltaasplinequ}) coincides with (\ref{ballfunctBSinequ}) which, thus, is a member of the family of inequalities of Theorem~\ref{thm:ineqVolProdFunc}.

\item[(b)] If $\mu$ is taken to be the uniform spherical propability measure $\sigma$, the induced Asplund endomorphism $\Psi_{\sigma}$ has the following interesting properties:
\begin{itemize}
\item $\Psi_{\sigma} f$ is radially symmetric for every $f \in \mathrm{LC_c}(\mathbb{R}^n)$;
\item $\Psi_{\sigma} \mathbbm{1}_K = \frac{w(K)}{2} \cdot \mathbbm{1}_{B^n}$ for every $K \in \mathcal{K}^n$.
\end{itemize}
Moreover, we will see that inequality (\ref{functinequAsplendos}) for $\Psi_{\sigma}$ is \emph{strictly} stronger than the functional analogue of Urysohn's inequality (\ref{functurysohn}) and, when restricted to indicators of convex bodies, yields a version of the Urysohn inequality (\ref{Uryinequ}) which is asymptotically sharp (see Section 4).
\end{enumerate}

\vspace{0.2cm}

The proof of Theorem \ref{thm:ineqVolProdFunc} does not make use of Theorem \ref{thm:ineqVolProdMonotone}, however, it follows similar arguments as in the geometric setting, replacing the application of the Blaschke--Santal\'o inequality (\ref{BSinequ}) by its functional form.
In particular, in Section 4 we also prove a functional version of Theorem~\ref{thm:ineqUrysohnMonotone} showing that each of the inequalities of Theorem \ref{thm:ineqVolProdFunc} is stronger than the one for the Asplund endomorphism $\Psi_{\sigma}$ and, hence, \emph{strictly} stronger than the functional analogue of Urysohn's inequality (\ref{functurysohn}). For even $\mu$, inequality (\ref{deltaasplinequ}) is proven to be the strongest one among the family of inequalities (\ref{functinequAsplendos}). Finally, we will see that Theorem \ref{thm:ineqVolProdMonotone} can be recovered in an \linebreak asymptotically optimal form by restricting Theorem \ref{thm:ineqVolProdFunc} to indicators of convex bodies.

\pagebreak

\centerline{\large{\bf{ \setcounter{abschnitt}{2}
\arabic{abschnitt}. Background material}}}

\reseteqn \alpheqn \setcounter{theorem}{0}

\vspace{0.6cm}

In this section we recall additional basic facts from convex geometry and some notions and results about convex and log-concave functions as well as the definition of the convolution of spherical functions and measures.
As general references, we recommend the monographs by Gardner \textbf{\cite{gardner2ed}}, Schneider \textbf{\cite{schneider93}}, Rockafellar~\textbf{\cite{Rockafellar1970}}, and Rockafellar and Wets \textbf{\cite{Rockafellar1998}}, as well as the survey \textbf{\cite{colesanti2017}} by Colesanti.

First recall that $\mathcal{K}^n$ denotes the space of convex bodies in $\mathbb{R}^n$ endowed with the Hausdorff metric. Each $K \in \mathcal{K}^n$ is uniquely determined by its support function $h(K,x) = \max \{x \cdot y: y \in K\}$, $x \in \mathbb{R}^n$, which is positively homogeneous of degree one and subadditive. Conversely, every function on $\mathbb{R}^n$ satisfying these two properties is the support function of a unique convex body.
For $K, L \in \mathcal{K}^n$, the support function of their Minkowski sum $K + L = \{x + y: x \in K, y \in L\}$ is given by
\begin{equation} \label{minkaddsuppfct}
h(K + L,\cdot) = h(K,\cdot) + h(L,\cdot).
\end{equation}
Moreover, for every $\vartheta \in \mathrm{SO}(n)$ and $y \in \mathbb{R}^n$, we have
\begin{equation} \label{sontranslsuppfct}
h(\vartheta K,x) = h(K,\vartheta^{-1}x) \qquad \mbox{and} \qquad h(K + y,x) = h(K,x) + x \cdot y
\end{equation}
for all $x \in \mathbb{R}^n$. Also the Hausdorff distance $d(K,L)$ of two convex bodies $K, L \in \mathcal{K}^n$ can be expressed conveniently by $d(K,L) = \|h(K,\,\cdot\,) - h(L,\,\cdot\,)\|_{\infty}$,
where $\|\,\cdot\,\|_{\infty}$ denotes the maximum norm on $C(\mathbb{S}^{n-1})$. Finally, recall that $K \subseteq L$ if and only if $h(K,\,\cdot\,) \leq h(L,\,\cdot\,)$, in particular, $h(K,\,\cdot\,) > 0$ if and only if $o \in \mathrm{int}\,K$.

Recall that if $K \in \mathcal{K}^n$ contains the origin in its interior, then $K^{\circ}$ denotes the polar body of $K$ with respect to the origin. We will make frequent use of the following polar coordinate formula for the volume of $K^{\circ}$,
\begin{equation} \label{polarformvol}
|K^{\circ}| = \frac{1}{n} \int_{\mathbb{S}^{n-1}} h(K,u)^{-n}\,du,
\end{equation}
where integration is with respect to spherical Lebesgue measure. The \emph{mean width} of a convex body $K \in \mathcal{K}^n$ is defined by
\begin{equation} \label{defmeanwidth}
w(K) = \frac{2}{n|B^n|}\int_{\mathbb{S}^{n-1}} h(K,u)\,du.
\end{equation}
The \emph{Steiner point} $s(K) \in \mathbb{R}^n$ of $K$ is the unique point in $\mathrm{relint}\,K$ defined by
\begin{equation} \label{defsteinerpkt}
s(K) = \frac{1}{|B^n|} \int_{\mathbb{S}^{n-1}} h(K,u)u\,du.
\end{equation}

Recall that the mean width and the Steiner point are uniquely determined by their Minkowski additivity and compatibility with rigid motions. To be more precise, a continuous map  $\varrho: \mathcal{K}^n \rightarrow \mathbb{R}$ is Minkowski additive and rigid motion invariant if and only if it is a constant multiple of the mean width $w$, while the Steiner point is the unique continuous map $s: \mathcal{K}^n \rightarrow \mathbb{R}^n$ which is Minkowski additive and rigid motion equivariant (cf.\ \textbf{\cite[\textnormal{Section 3.3}]{schneider93}}).

\pagebreak

Associated with each convex body $K \in \mathcal{K}^n$ is a finite non-negative Borel measure $S_1(K,\cdot)$ on $\mathbb{S}^{n-1}$, its \emph{area measure of order one}. Using the Laplacian (or Laplace--Beltrami operator) $\Delta_{\mathbb{S}}$ on $\mathbb{S}^{n-1}$, it can be defined by the relation (to be understood in a distributional sense)
\[S_1(K,\cdot) = \frac{1}{n-1}\, \Delta_{\mathbb{S}}\,h(K,\cdot) + h(K,\cdot).  \]
As shown by Weil \textbf{\cite{weil80}}, the set $\{S_1(K,\cdot): K \in \mathcal{K}^n\}$ is \emph{not} dense in the set of all non-negative measures on $\mathbb{S}^{n-1}$ with barycenter at the origin in the weak topology.

Each non-negative even measure $\mu$ on $\mathbb{S}^{n-1}$ generates a uniquely determined origin-symmetric convex body $Z^{\mu} \in \mathcal{K}^n$ by
\begin{equation} \label{zonoidmeas}
h(Z^{\mu},u) = \int_{\mathbb{S}^{n-1}} |u \cdot v|\,d\mu(v), \qquad u \in \mathbb{S}^{n-1}.
\end{equation}
The bodies obtained in this way constitute the class of origin-symmetric \emph{zonoids}, which naturally arise also in various other contexts (see, e.g., \textbf{\cite[\textnormal{Chapter~3.5}]{schneider93}}).

\vspace{0.2cm}

We turn now to convex and log-concave functions on $\mathbb{R}^n$. Let $\mathrm{Cvx}(\mathbb{R}^n)$ denote the set of convex and lower semi-continuous functions $\varphi: \mathbb{R}^n \rightarrow (-\infty,\infty]$ which are \emph{proper}, that is, not identically $+\infty$. Two convex sets naturally associated to any $\varphi \in \mathrm{Cvx}(\mathbb{R}^n)$ are its \emph{domain},
$\mathrm{dom}\,\varphi = \{x \in \mathbb{R}^n: \varphi(x) < +\infty\}$,
and its \emph{epigraph} defined by $\mathrm{epi}\,\varphi = \{(x,\xi) \in \mathbb{R}^n \times \mathbb{R}: \varphi(x) \leq \xi\}$.
Note that for any $\varphi \in \mathrm{Cvx}(\mathbb{R}^n)$, $\mathrm{dom}\,\varphi$ is non-empty and $\mathrm{epi}\,\varphi$ is closed and non-empty.
We call a function $\varphi \in \mathrm{Cvx}(\mathbb{R}^n)$ \emph{coercive} if $\lim_{\|x\|\rightarrow \infty} \varphi(x) = +\infty$ and we denote by $\mathrm{Cvx_c}(\mathbb{R}^n)$
the set of all coercive $\varphi \in \mathrm{Cvx}(\mathbb{R}^n)$. We also note that (see, e.g., \textbf{\cite[\textnormal{Lemma 2.5}]{ColFrag2013}}),
$\varphi \in \mathrm{Cvx}(\mathbb{R}^n)$ is coercive if and only if there exist $\gamma > 0$ and $\beta \in \mathbb{R}$ such that for every $x \in \mathbb{R}^n$,
\begin{equation} \label{coercconecond}
\varphi(x) \geq \gamma \|x\| + \beta.
\end{equation}
Next we endow the spaces $\mathrm{Cvx}(\mathbb{R}^n)$ and $\mathrm{Cvx_c}(\mathbb{R}^n)$ with the topology induced by epi-convergence. Recall that a sequence of $\varphi_k \in \mathrm{Cvx}(\mathbb{R}^n)$ is called
\emph{epi-convergent} to $\varphi: \mathbb{R}^n \rightarrow (-\infty,\infty]$ if for all $x \in \mathbb{R}^n$ the following two conditions hold:
\begin{itemize}
\item $\varphi(x) \leq {\liminf}_{k \rightarrow \infty} \varphi_k(x_k)$ for every sequence $x_k$ that converges to $x$.
\item There exists a sequence $x_k$ converging to $x$ such that $\varphi(x) = \lim_{k \rightarrow \infty} \varphi_k(x_k)$.
\end{itemize}
In this case, we write $\varphi_k \stackrel{\mathrm{epi}}{\rightarrow} \varphi$. Note that the limiting function $\varphi$ is again convex and lower semi-continuous.
However, in general, $\varphi$ need not be proper.

\begin{lem} \label{lemepiconvequiv} \emph{(\!\!\textbf{\cite[\textnormal{Theorem 7.17}]{Rockafellar1998}})} If $\varphi, \varphi_k \in \mathrm{Cvx}(\mathbb{R}^n)$ and $\mathrm{int}\,\mathrm{dom}\,\varphi$ is non-empty, then the following statements are equivalent to $\varphi_k$ being epi-convergent to~$\varphi$:
\begin{enumerate}
\item[(i)] There exists a dense set $D \subseteq \mathbb{R}^n$ such that $\varphi_k(x) \rightarrow \varphi(x)$ for every $x \in D$.
\item[(ii)] The sequence $\varphi_k$ converges uniformly to $\varphi$ on every compact subset of $\mathbb{R}^n$ that does not intersect the boundary of $\mathrm{dom}\, \varphi$.
\end{enumerate}
\end{lem}

\pagebreak

\noindent Let us also emphasize that epi-convergence (also known as $\Gamma$-convergence) is equivalent to the convergence of the corresponding epigraphs in the socalled Painlev\'e--Kuratowski sense
(cf.\ \textbf{\cite[\textnormal{Section 7.B}]{Rockafellar1998}}).

For $\varphi, \psi \in \mathrm{Cvx}(\mathbb{R}^n)$, their \emph{infimal convolution} is defined by
\[ (\varphi\, \mbox{{\scriptsize $\square$}}\, \psi)(x) = \inf_{x_1 + x_2 = x} \{\varphi(x_1) + \psi(x_2)\}.  \]
If $\varphi\, \mbox{{\scriptsize $\square$}}\, \psi$ does not attain $-\infty$, then it is convex, proper, and $\mathrm{epi}(\varphi\, \mbox{{\scriptsize $\square$}}\, \psi) = \mathrm{epi}\,\varphi + \mathrm{epi}\,\psi$. 
However, $\varphi\, \mbox{{\scriptsize $\square$}}\, \psi$ need not be semi-continuous. A quite useful condition to ensure lower semi-continuity of the infimal
convolution of $\varphi, \psi \in \mathrm{Cvx}(\mathbb{R}^n)$ can be found in \textbf{\cite[\textnormal{Corollary 9.2.2}]{Rockafellar1970}} and requires that
\begin{equation} \label{semicontcondit}
\lim_{\lambda \rightarrow \infty} \frac{\varphi(y + \lambda x)}{\lambda} + \lim_{\lambda \rightarrow \infty} \frac{\psi(z - \lambda x)}{\lambda}>0
\end{equation}
for every non-zero $x \in \mathbb{R}^n$ and arbitrary $y \in \mathrm{dom}\,\varphi$, $z \in \mathrm{dom}\,\psi$.

For $t > 0$ and $\varphi \in \mathrm{Cvx}(\mathbb{R}^n)$, the \emph{Moreau envelope} $e_t \varphi$ of $\varphi$ is defined by
\[ e_t \varphi = \varphi \, \mbox{{\scriptsize $\square$}}\, \mbox{$\frac{1}{2t}$} \|\cdot\|^2.\]
For the proof of Theorem \ref{def:phimuFirst}, we require the following of its simple properties.

\begin{lem}\label{propMoreau} \emph{(\!\!\textbf{\cite[\textnormal{Theorems 1.25 \& 2.26}]{Rockafellar1998}})}
Suppose that $\varphi \in \mathrm{Cvx}(\mathbb{R}^n)$. Then the following statements hold:
 \begin{itemize}
  \item[(i)] \!\! $e_t \varphi \in \mathrm{Cvx}(\mathbb{R}^n)$ and it is finite for every $t > 0$;
  \item[(ii)] \!\! $e_t \varphi(x)$ converges to $\varphi(x)$ monotonously from below for every $x \in \mathbb{R}^n$  as $t \searrow 0$.
 \end{itemize}
In particular, $e_t \varphi \stackrel{\mathrm{epi}}{\rightarrow} \varphi$ as $t \searrow 0$.
\end{lem}

Now, let $\mathrm{LC}(\mathbb{R}^n) = \{f = e^{-\varphi}: \varphi \in \mathrm{Cvx}(\mathbb{R}^n)\}$ denote the set of all proper, log-concave, and upper semi-continuous functions on
$\mathbb{R}^n$ and recall that
\[\mathrm{LC_c}(\mathbb{R}^n) = \left \{f \in \mathrm{LC}(\mathbb{R}^n): \lim_{\|x\|\rightarrow \infty} f(x) = 0\right\} =
\left \{f = e^{-\varphi}: \varphi \in \mathrm{Cvx_c}(\mathbb{R}^n)\right\}.\]
We call a sequence $f_k = e^{-\varphi_k} \in \mathrm{LC}(\mathbb{R}^n)$ or $\mathrm{LC_c}(\mathbb{R}^n)$, respectively, hypo-convergent to $f = e^{-\varphi}$,
if $\varphi_k \in \mathrm{Cvx}(\mathbb{R}^n)$ or $\mathrm{Cvx_c}(\mathbb{R}^n)$, respectively, epi-converges to $\varphi$. (Note that for log-concave functions our notion of hypo-convergence coincides with the more general definition used frequently in analysis.)

Next, recall that the Asplund sum $f \star g$ of $f = e^{-\varphi}, g = e^{-\psi} \in \mathrm{LC}(\mathbb{R}^n)$ is related to the infimal convolution $\varphi\, \mbox{{\scriptsize $\square$}}\, \psi$ of
$\varphi, \psi \in \mathrm{Cvx}(\mathbb{R}^n)$ by
\begin{equation} \label{asplinfconv}
f \star g = e^{- \varphi\, \mbox{{\scriptsize $\square$}}\, \psi}.
\end{equation}
In particular, since $\mathrm{Cvx}(\mathbb{R}^n)$ is not closed under infimal convolution, the space $\mathrm{LC}(\mathbb{R}^n)$ is not closed under Asplund addition.
However, as the following lemma shows, this is no longer the case when considering coercive functions.

\begin{lem} \label{asplproplem}
Suppose that $f, g \in \mathrm{LC_c}(\mathbb{R}^n)$, $a, b > 0$ and $y \in \mathbb{R}^n$. Then the following statements hold:
 \begin{itemize}
  \item[(i)] $f \star g \in \mathrm{LC_c}(\mathbb{R}^n)$;
  \item[(ii)] $(a \cdot f) \star (b \cdot f) = (a + b) \cdot f$;
  \item[(iii)] $(f \star \mathbbm{1}_{\{y\}})(x) = f(x - y)$.
 \end{itemize}
 \end{lem}

\noindent {\it Proof.} In order to prove (i), let $f = e^{-\varphi}$ and $g = e^{-\psi}$ with $\varphi, \psi \in \mathrm{Cvx_c}(\mathbb{R}^n)$. Then, by (\ref{asplinfconv}),
it is sufficient to show that $\varphi\, \mbox{{\scriptsize $\square$}}\, \psi \in \mathrm{Cvx_c}(\mathbb{R}^n)$. First note that since $\varphi$ and $\psi$ are coercive, they are
bounded from below and, consequently, so is $\varphi\, \mbox{{\scriptsize $\square$}}\, \psi$ which is, therefore, convex and proper. Moreover, from the definition of infimal convolution and the triangle
inequality, it follows that $\varphi\, \mbox{{\scriptsize $\square$}}\, \psi$ is coercive. It remains to show that $\varphi\, \mbox{{\scriptsize $\square$}}\, \psi$ is
lower semi-continuous, for which we check that condition (\ref{semicontcondit}) is satisfied. To this end, we use (\ref{coercconecond}) to conclude that there exist
$\gamma_{\varphi}, \gamma_{\psi} > 0$ and $\beta_{\varphi}, \beta_{\psi} \in \mathbb{R}$ such that $\varphi(w) \geq \gamma_{\varphi}\|w\| + \beta_{\varphi}$ and
$\psi(w) \geq \gamma_{\psi}\|w\| + \beta_{\psi}$ for every $w \in \mathbb{R}^n$. Thus,
\[\lim_{\lambda \rightarrow \infty} \frac{\varphi(y + \lambda x)}{\lambda} + \lim_{\lambda \rightarrow \infty} \frac{\psi(z - \lambda x)}{\lambda}
\geq \gamma_{\varphi}\|x\| + \gamma_{\psi}\|x\| > 0\]
for every non-zero $x \in \mathbb{R}^n$ which completes the proof of (i). Statements (ii) and (iii) follow easily from the definition of the Asplund sum and multiplication.
\hfill $\blacksquare$

\vspace{0.5cm}

\noindent {\bf Examples:}

\vspace{-0.1cm}

\begin{enumerate}
\item[(a)] If $K \in \mathcal{K}^n$, then the indicator function $\mathbbm{1}_K$ belongs to $\mathrm{LC_c}(\mathbb{R}^n)$.

\item[(b)] Recall that for $K \in \mathcal{K}^n$ containing the origin in its interior, its \emph{gauge} or \emph{Minkowski functional} is given by
\[\|x\|_K=\min\{\lambda \geq 0: x \in \lambda K\}, \qquad x \in \mathbb{R}^n.  \]
If $K$ is origin-symmetric, $\|\cdot\|_K$ is the norm with unit ball $K$. For $K = B^n$, we have $\|\cdot\|_{B^n} = \|\cdot\|$.
Another interesting class of log-concave functions consists of those $f = e^{-\varphi} \in \mathrm{LC_c}(\mathbb{R}^n)$, where
\begin{equation} \label{normhochp}
\varphi = \frac{1}{p} \|\cdot\|_K^p, \qquad p \geq 1.
\end{equation}
In particular, $f \in \mathrm{LC_c}(\mathbb{R}^n)$ is called a \emph{Gaussian} if there exist $a > 0$, $y \in \mathbb{R}^n$ and an origin-symmetric ellipsoid $E \subseteq \mathbb{R}^n$ such that
\[f(x) = a\, e^{-\frac{1}{2}\|x-y\|_E^2}, \qquad x \in \mathbb{R}^n.  \]
For $a = (2\pi)^{-n/2}$, $y = o$ and $E = B^n$, we obtain the \emph{standard Gaussian} $\psi_n$.
\end{enumerate}

Let us mention here another useful volume formula for convex bodies, involving the functions defined in (\ref{normhochp}). If $K \in \mathcal{K}^n$ contains the origin in its interior, then
\begin{equation} \label{volform17}
|K| = \frac{1}{p^{n/p}\Gamma\left (1 + \frac{n}{p}\right )} \int_{\mathbb{R}^n} \exp\!\left (\! -\frac{1}{p} \|x\|_K^p  \right )dx.
\end{equation}

The \emph{Legendre transform} $\mathcal{L}: \mathrm{Cvx}(\mathbb{R}^n) \rightarrow \mathrm{Cvx}(\mathbb{R}^n)$ is defined by
\[(\mathcal{L}\varphi)(x) = \sup_{y \in \mathbb{R}^n} x \cdot y - \varphi(y), \qquad x \in \mathbb{R}^n.  \]
It is a classical notion with many applications in several areas which are extensively covered in the literature (e.g., \textbf{\cite{Rockafellar1970, Rockafellar1998}}).
We collect a number of its well known properties for quick later reference in the following proposition. Note that the properties from (i) were recently shown to
essentially characterize the Legendre transform in a fundamental paper by Artstein-Avidan and Milman \textbf{\cite{Artstein2009}}.

\begin{prop} \label{propLegendre} \hspace{-0.2cm} \emph{(see, e.g., \textbf{\cite{Rockafellar1970, Rockafellar1998}})} For $\varphi_k, \varphi, \psi \in \mathrm{Cvx}(\mathbb{R}^n)$ and $K \in \mathcal{K}^n$, the following statements hold:
\begin{itemize}
\item[(i)] $\mathcal{L}\mathcal{L}\varphi = \varphi$ and if $\varphi \leq \psi$, then $\mathcal{L} \varphi \geq \mathcal{L} \psi$;
\item[(ii)] $\varphi$ is coercive if and only if $\mathrm{dom}\,\mathcal{L} \varphi$ contains the origin in its interior;\\[-0.5cm]
\item[(iii)] $\varphi_k \stackrel{\mathrm{epi}}{\rightarrow} \varphi$ if and only if $\mathcal{L}\varphi_k \stackrel{\mathrm{epi}}{\rightarrow} \mathcal{L}\varphi$;
\item[(iv)] $\mathcal{L}(-\log \mathbbm{1}_K) = h(K,\cdot)$ and if $1 < p,q < \infty$ are such that $\frac{1}{p}+\frac{1}{q}=1$ and $K$ contains the origin in its interior, then
\[\mathcal{L}\left (\frac{1}{p} \|\cdot\|_K^p\right ) =  \frac{1}{q} \|\cdot\|_{K^\circ}^q. \]
\end{itemize}

\end{prop}

The Legendre transform gives rise to several constructions for log-concave functions. Recall that for $f \in \mathrm{LC}(\mathbb{R}^n)$, its support function is
$h(f,\cdot) = \mathcal{L}(-\log f)$ and note that $h(f,\cdot) \in \mathrm{Cvx}(\mathbb{R}^n)$. For our purposes it is particularly useful to observe that
support functions are Asplund additive (cf.\ \textbf{\cite[\textnormal{p.518}]{schneider93}}), in the sense that for $f, g \in \mathrm{LC_c}(\mathbb{R}^n)$, we have
\begin{equation} \label{suppfctaspladd}
h(f \star\, g,\,\cdot\,) = h(f,\,\cdot\,) + h(g,\,\cdot\,)
\end{equation}
which is an extension of (\ref{minkaddsuppfct}) to $\mathrm{LC_c}(\mathbb{R}^n)$. Moreover, for $\vartheta \in \mathrm{SO}(n)$ and $y \in \mathbb{R}^n$,
\begin{equation} \label{sontranslsuppfctlcc}
h(\vartheta f,x) = h(f,\vartheta^{-1}x) \qquad \mbox{and} \qquad h(f(\,.-y),x) = h(f,x) + x \cdot y
\end{equation}
for every $x \in \mathbb{R}^n$ -- an extension of properties (\ref{sontranslsuppfct}) to $\mathrm{LC_c}(\mathbb{R}^n)$.

Recall that for $f \in \mathrm{LC}(\mathbb{R}^n)$, its polar function is defined by $f^\circ = e^{-\mathcal{L}(-\log f)}$.
From the definition of the Legendre transform and Proposition \ref{propLegendre}, one obtains:

\pagebreak

\begin{lem} For $f,g \in \mathrm{LC_c}(\mathbb{R}^n)$, $a > 0$, and $K \in \mathcal{K}^n$ containing the origin in its interior, the following statements hold:
\begin{itemize}
\item[(i)] $(f^{\circ})^{\circ} = f$;
\item[(ii)] $(Af)^{\circ} = A^{-\mathrm{T}}f^{\circ}$ for every $A \in \mathrm{GL}(n)$;
\item[(iii)] $(f \star g)^{\circ} = f^\circ\,g^\circ$ and $(a \cdot f)^{\circ} = (f^\circ)^a$;
\item[(iv)] $(\mathbbm{1}_K)^{\circ} = e^{-\|\cdot\|_{K^{\circ}}}$ and if $1 < p,q < \infty$ are such that $\frac{1}{p}+\frac{1}{q}=1$, then
\[ \left ( e^{-\frac{1}{p} \|\cdot\|_K^p}  \right )^{\circ} = e^{-\frac{1}{q}\|\cdot\|_{K^{\circ}}^q}. \]
\end{itemize}
\end{lem}

We next state a more general version of the functional Blaschke--Santal\'o inequality due to Artstein-Avidan, Klartag, and Milman \textbf{\cite{Artstein2004}} required in the proof of Theorem \ref{thm:ineqVolProdFunc}.
We restrict ourselves to log-concave functions and recall that the centroid of an integrable function $f$ on $\mathbb{R}^n$ such that
$\int_{\mathbb{R}^n}f\,dx>0$ is defined by
\[\mathrm{cent}\,f = \frac{\int_{\mathbb{R}^n}x f(x)\,dx}{\int_{\mathbb{R}^n} f(x)\,dx}.  \]

\begin{theorem} \label{generalfctBS} \emph{(Artstein-Avidan et al.\ \textbf{\cite{Artstein2004}})}  Suppose that $f \in \mathrm{LC}(\mathbb{R}^n)$ is such that $0 < \int_{\mathbb{R}^n} f(x)\,dx < \infty$ and let $\tilde{f}(x) = f(x - \mathrm{cent}\,f)$. Then
\[\int_{\mathbb{R}^n} f(x)\,dx  \int_{\mathbb{R}^n} \tilde{f}^{\circ}(x)\,dx \leq (2\pi)^n \]
with equality if and only if $f$ is a Gaussian.
\end{theorem}

Let us note that for $K \in \mathcal{K}^n$ containing the origin in its interior and $f = e^{-\frac{1}{2}\|\cdot\|_K^2}$ in Theorem \ref{generalfctBS}, one
recovers a version of the geometric Blaschke--Santal\'o inequality equivalent to (\ref{BSinequ}).

\vspace{0.2cm}

In the final part of this section, we recall the convolution of spherical functions and measures and its extension to functions on $\mathbb{R}^n$.
We are particularly interested in convolutions with \emph{zonal} measures, that is, $\mathrm{SO}(n-1)$ invariant measures on $\mathbb{S}^{n-1}$, where $\mathrm{SO}(n-1)$ is the subgroup of $\mathrm{SO}(n)$ keeping an arbitrary pole $\bar{e} \in \mathbb{S}^{n-1}$ fixed. First, recall that since $\mathrm{SO}(n)$ is a compact Lie group, the convolution $\tau \ast \mu$ of signed measures $\tau, \mu$ on $\mathrm{SO}(n)$ can be defined by
\[\int_{\mathrm{SO}(n)}\!\!\! f(\vartheta)\, d(\tau \ast \mu)(\vartheta)=\int_{\mathrm{SO}(n)}\! \int_{\mathrm{SO}(n)}\!\!\! f(\eta \theta)\,d\tau(\eta)\,d\mu(\theta), \qquad f \in C(\mathrm{SO}(n)).   \]
Since $\mathbb{S}^{n-1}$ is diffeomorphic to the homogeneous space $\mathrm{SO}(n)/\mathrm{SO}(n-1)$, there is a natural one-to-one correspondence between functions and measures on $\mathbb{S}^{n-1}$ and right $\mathrm{SO}(n-1)$ invariant functions and measures on $\mathrm{SO}(n)$, respectively. Using this correspondence, the convolution of measures on $\mathrm{SO}(n)$ induces an associative convolution product of measures on $\mathbb{S}^{n-1}$ (cf.\ \textbf{\cite{Schu09}} for more details).

\pagebreak

Here, we only note that the spherical convolution with a zonal measure on $\mathbb{S}^{n-1}$ takes an especially simple form. If for $u \in \mathbb{S}^{n-1}$, we denote by $\vartheta_u \in \mathrm{SO}(n)$ an arbitrary rotation such that $ \vartheta_u\bar{e} = u$, the convolution of $h \in C(\mathbb{S}^{n-1})$ and a zonal measure $\mu$ on $\mathbb{S}^{n-1}$ is given by
\begin{equation} \label{zonalconv}
(h \ast \mu)(u) = \int_{\mathbb{S}^{n-1}} h(\vartheta_u v)\,d\mu(v), \qquad u \in \mathbb{S}^{n-1}.
\end{equation}
Since $\vartheta_{\eta u} = \eta \vartheta_u$ for any $\eta \in \mathrm{SO}(n)$ up to right-multiplication by an element of $\mathrm{SO}(n-1)$, it follows from (\ref{zonalconv}) that
$(\eta h) \ast \mu = \eta(h \ast \mu)$, where $(\eta h)(u) = h(\eta^{-1}u)$, $u \in \mathbb{S}^{n-1}$. It is also not difficult to check that the convolution (\ref{zonalconv}) is selfadjoint and that spherical convolution of zonal functions and measures is Abelian.

Motivated by (\ref{zonalconv}) and the significance of the spherical convolution  $\ast$ for Minkowski endomorphisms, we now introduce an extension of $\ast$ to the following important open subset of convex functions in $\mathrm{Cvx}(\mathbb{R}^n)$,
\[\mathrm{Cvx}_{(o)}(\mathbb{R}^n) = \{\varphi \in \mathrm{Cvx}(\mathbb{R}^n): o \in \mathrm{int}\, \mathrm{dom}\,\varphi \}.  \]
Note that, by Proposition \ref{propLegendre} (ii), $\varphi \in \mathrm{Cvx}_{(o)}$ if and only if $\mathcal{L}\varphi \in \mathrm{Cvx_c}(\mathbb{R}^n)$
or, equivalently,
\begin{equation} \label{lccundcvx0}
f \in \mathrm{LC_c}(\mathbb{R}^n) \mbox{ if and only if } h(f,\cdot) \in \mathrm{Cvx}_{(o)}(\mathbb{R}^n).
\end{equation}
In the following, for $x \in \mathbb{R}^n\backslash \{o\}$, let $\vartheta_x \in \mathrm{SO}(n)$ denote an arbitrary rotation such that $\vartheta_x \bar{e} = \frac{x}{\|x\|}$.

\vspace{0.3cm}

\noindent {\bf Definition.} Suppose that $\varphi \in \mathrm{Cvx}_{(o)}(\mathbb{R}^n)$ and let $\mu$ be a non-negative zonal measure on $\mathbb{S}^{n-1}$. The convolution $\varphi \circledast  \mu$ is defined for $x \in \mathbb{R}^n\backslash \{o\}$ by
\begin{equation} \label{defcircledast}
  (\varphi \circledast \mu) (x) = \int_{\mathbb{S}^{n-1}} \varphi(\|x\| \vartheta_x v) \, d\mu(v)
\end{equation}
and at the origin by $(\varphi \circledast \mu)(o) = \liminf_{\|x\| \rightarrow 0}(\varphi \circledast \mu) (x)$ .

\vspace{0.3cm}

We will show in the next section that $\varphi \circledast \mu$ is indeed a well defined function in $\mathrm{Cvx}_{(o)}(\mathbb{R}^n)$, but already note here that if $\varphi$ is homogeneous of some degree $p \in \mathbb{R}$, then, by (\ref{zonalconv}), $\varphi \circledast \mu$ coincides with the homogeneous extension of degree $p$ of $\widehat{\varphi} \ast \mu$ to $\mathbb{R}^n$, where $\widehat{\varphi}$ is the restriction of $\varphi$ to $\mathbb{S}^{n-1}$.

\vspace{1cm}

\centerline{\large{\bf{ \setcounter{abschnitt}{3}
\arabic{abschnitt}. Minkowski and Asplund endomorphisms}}}

\reseteqn \alpheqn \setcounter{theorem}{0}

\vspace{0.6cm}

In the following we first review additional background material on Minkowski endomorphisms required in the proofs of Theorems \ref{thm:ineqVolProdMonotone} and
\ref{thm:ineqUrysohnMonotone} and the discussion thereof in Section 4. In the second part of this section we give the proof of Theorem~\ref{def:phimuFirst} as well as that of an auxiliary result
used in the proof of an analogue of Theorem~\ref{thm:ineqUrysohnMonotone} for log-concave functions.

\pagebreak

We begin by recalling a notion of monotonicity for operators on convex bodies which is of particular importance for Minkowski endomorphisms.

\vspace{0.3cm}

\noindent {\bf Definition} A map $\Phi: \mathcal{K}^n \rightarrow \mathcal{K}^n$ is called \emph{weakly monotone} if $\Phi K \subseteq \Phi L$ for all $K, L \in \mathcal{K}^n$ such
that $K \subseteq L$ and $s(K) = s(L) = o$.

\vspace{0.3cm}

The quest to establish a classification of all Minkowski endomorphisms has its origin in the paper \textbf{\cite{Schneider1974c}} from 1974 by Schneider. The following result -- combining theorems by Dorrek and Kiderlen -- represents the status quo on this difficult task which has not yet been completed. Here, we call a measure on $\mathbb{S}^{n-1}$ \emph{linear} if it has a density (w.r.t.\ spherical Lebesgue measure) of the form $u \mapsto x \cdot u$ for some $x \in \mathbb{R}^n$.

\begin{theorem} \label{weakmonMinkendoclass} \emph{(Dorrek \textbf{\cite{Dorrek2017b}}, Kiderlen \textbf{\cite{Kiderlen2006}})} If $\Phi: \mathcal{K}^n \rightarrow \mathcal{K}^n$ is a Minkowski
endomorphism, then there exists a signed zonal measure $\mu$ on $\mathbb{S}^{n-1}$ with center of mass at the origin such that
\begin{equation} \label{faltdarstell17}
h(\Phi K,\cdot) = h(K,\cdot) \ast \mu
\end{equation}
for every $K \in \mathcal{K}^n$. The measure $\mu$ is uniquely determined by $\Phi$.\\[0.1cm] Moreover, $\Phi$ is monotone if and only if $\mu$ is non-negative and $\Phi$ is weakly
monotone if and only if $\mu$ is non-negative up to addition of a linear measure.
\end{theorem}

The measure $\mu$ uniquely associated with the Minkowski endomorphism $\Phi$ via the relation (\ref{faltdarstell17}) is called the \emph{generating measure} of $\Phi$ and we frequently indicate this by writing $\Phi_{\mu}$. Exploiting this one-to-one correspondence, we can endow the cone of Minkowski endomorphisms with the topology induced by the weak convergence of their generating measures. Before we exhibit some prominent examples, let us note that for $n = 2$, Theorem \ref{weakmonMinkendoclass} is due to Schneider \textbf{\cite{Schneider1974}} who also showed that in this special case \emph{all} Minkowski endomorphisms are weakly monotone. The conjecture that the same is true for $n \geq 3$ was disproved by Dorrek \textbf{\cite{Dorrek2017b}}.

\vspace{0.3cm}

\noindent {\bf Examples:}

\begin{enumerate}
\item[(a)] Recall that $\sigma$ denotes the $\mathrm{SO}(n)$ invariant probability measure on $\mathbb{S}^{n-1}$.
The Minkowski endomorphism $\Phi_{\sigma}: \mathcal{K}^n \rightarrow \mathcal{K}^n$ generated by $\sigma$ satisfies
\begin{equation} \label{meanwidthendo}
\Phi_{\sigma}K = \frac{w(K)}{2} B^n
\end{equation}
for every $K \in \mathcal{K}^n$ and, thus, the inequality $|K||\Phi_{\sigma}^\circ K| \leq |B^n|^2$ is precisely the Urysohn inequality.

\item[(b)] The unique \emph{discrete} zonal probability measure on $\mathbb{S}^{n-1}$ with center of mass at the origin is given by
\begin{equation} \label{discrzonal}
\nu = \frac{1}{2}(\delta_{\bar{e}} + \delta_{-\bar{e}}).
\end{equation}
It is the generating measure of the central symmetrization $\Delta: \mathcal{K}^n \rightarrow \mathcal{K}^n$.

\item[(c)] The generating measure of the Minkowski endomorphism $\Pi_1: \mathcal{K}^n \rightarrow \mathcal{K}^n$ (recall our normalization $\Pi_1 B^n = B^n$)
is the invariant probability measure $\sigma_{\bar{e}^\perp}$ concentrated on the equator $\mathbb{S}^{n-1} \cap \bar{e}^\perp$. Noting that
\[\sigma_{\bar{e}^\perp} = S_1\!\left (\mbox{$\frac{1}{2}$}[-\bar{e},\bar{e}],\cdot\right)\!,    \]
Berg and the second author \textbf{\cite{bergschu2020}} considered, more generally, Minkowski endomorphisms generated by measures of the form $S_1(Z,\cdot)$ for some zonoid $Z \in \mathcal{K}^n$ and established Theorem \ref{thm:ineqVolProdMonotone} for such maps. However, this class is not dense in all monotone Minkowski endomorphisms.

\item[(d)] The Minkowski endomorphism $\mathrm{J}: \mathcal{K}^n \rightarrow \mathcal{K}^n$, defined by
\begin{equation} \label{idendo}
\mathrm{J}K = K - s(K),
\end{equation}
is weakly monotone and its generating measure is given by $\delta_{\bar{e}} - n(\bar{e} \cdot.\,)d\sigma$.
\end{enumerate}

For the reader's convenience, we prove next two useful previously known properties of Minkowski endomorphisms required in the proof of Theorem \ref{thm:ineqUrysohnMonotone}.

\begin{lem} \label{meanwidthminkendo} Suppose that $\Phi: \mathcal{K}^n \rightarrow \mathcal{K}^n$ is a Minkowski endomorphism with generating measure $\mu$ on $\mathbb{S}^{n-1}$.
Then the following statements hold:
\begin{enumerate}
\item[(i)] $w(\Phi K) = \mu(\mathbb{S}^{n-1})\, w(K)$ for every $K \in \mathcal{K}^n$;
\item[(ii)] If $\Phi$ is non-trivial and weakly monotone and $K \in \mathcal{K}^n$ has non-empty interior, then $\Phi K$ contains the origin in its interior;
\item[(iii)] $\Phi\{x\} = \{o\}$ for every $x \in \mathbb{R}^n$.
\end{enumerate}
\end{lem}

\noindent {\it Proof.} In order to see (i), we use that the spherical convolution is selfadjoint and Abelian for zonal measures. These facts combined with (\ref{defmeanwidth}) and (\ref{meanwidthendo}) yield,
\begin{align*}
w(\Phi K) & =  2 \int_{\mathbb{S}^{n-1}} (h(K,\cdot)\ast \mu)(u) \,d\sigma(u) = 2 \int_{\mathbb{S}^{n-1}} (h(K,\cdot)\ast \sigma)(u) \,d\mu(u) \\
& = w(K) \int_{\mathbb{S}^{n-1}} h(B^n,u)\,d\mu(u) = \mu(\mathbb{S}^{n-1})\, w(K)
\end{align*}
for every $K \in \mathcal{K}^n$. For the proof of (ii), first note that every non-trivial Minkowski endomorphism maps Euclidean balls of positive radii to origin-symmetric balls by $\mathrm{SO}(n)$ equivariance and translation invariance. These balls must be of positive radii by (i).
Now, using that the Steiner point $s(K) \in \mathrm{int}\,K$ for every $K \in \mathcal{K}^n$ with non-empty interior (see, e.g., \textbf{\cite[\textnormal{p.50}]{schneider93}}), 
we obtain (ii) from the translation invariance of $\Phi$, which implies $\Phi K = \Phi(K - s(K))$, and its monotonicity on bodies with Steiner points at the origin. Claim (iii) is a direct consequence of the fact that $\{x\} = \{x\} + \{o\}$ and the additivity and translation invariance of $\Phi$.
\hfill $\blacksquare$

\vspace{0.3cm}

We turn now to Asplund endomorphisms and prove the following critical result at the core of Theorem \ref{def:phimuFirst}.

\begin{prop} \label{propforthm3} Suppose that $\varphi \in \mathrm{Cvx}_{(o)}(\mathbb{R}^n)$ and $\mu$ is a non-negative zonal measure on $\mathbb{S}^{n-1}$. Then
$\varphi \circledast \mu$ is a well defined function in $\mathrm{Cvx}_{(o)}(\mathbb{R}^n)$.
Moreover, the map $\varphi \mapsto \varphi \circledast \mu$ defines a continuous linear operator from $\mathrm{Cvx}_{(o)}(\mathbb{R}^n)$ to itself which commutes with the action of
$\mathrm{SO}(n)$.
\end{prop}

\noindent {\it Proof.} First note that the right hand side of (\ref{defcircledast}) is independent of the choice of $\vartheta_x$ by the $\mathrm{SO}(n-1)$ invariance of $\mu$. Since $\varphi$ is convex on $\mathbb{R}^n$, it is bounded from below by an affine function and, thus, the negative part of $\varphi$ is bounded on every sphere in $\mathbb{R}^n$.
Consequently, the integral in (\ref{defcircledast}) is well defined and takes values in $(-\infty,\infty]$. Moreover, since $o \in \mathrm{int}\,\mathrm{dom}\,\varphi$, there exists $r > 0$ such that $\varphi$ takes finite values on $rB^n$, which implies that also $\varphi \circledast \mu$ takes finite values on $rB^n$. In particular, $\varphi \circledast \mu$ is proper and $o \in \mathrm{int}\,\mathrm{dom}\,\varphi \circledast \mu$.

The proof that $\varphi \circledast \mu$ is convex on $\mathbb{R}^n$ is rather tedious and technical and we therefore postpone it to the Appendix. In order to see that $\varphi \circledast \mu \in \mathrm{Cvx}_{(o)}(\mathbb{R}^n)$, it remains to show that it is lower semi-continuous. To this end, let $x_0 \in \mathbb{R}^n\backslash\{o\}$ and note that $\varphi$ is bounded from below on the compact set $2\|x_0\|B^n$ by semi-continuity. Hence, we may use Fatou's Lemma and the semi-continuity of $\varphi$ to conclude that
\begin{align*}
\liminf_{x \rightarrow x_0} \,(\varphi \circledast \mu)(x) & \geq \int_{\mathbb{S}^{n-1}} \liminf_{x \rightarrow x_0} \varphi(\|x\|\vartheta_x u)\,d\mu(u) \\
& \geq \int_{\mathbb{S}^{n-1}} \varphi(\|x_0\|\vartheta_{x_0} u)\,d\mu(u) = (\varphi \circledast \mu)(x_0)
\end{align*}
which combined with our definition of $(\varphi \circledast \mu)(o)$ yields the lower semi-continuity of $\varphi \circledast \mu$ on $\mathbb{R}^n$ and completes the proof that $\varphi \circledast \mu \in \mathrm{Cvx}_{(o)}(\mathbb{R}^n)$.

Since the linearity and commutativity with respect to the action of $\mathrm{SO}(n)$ of the map $\varphi \mapsto \varphi \circledast \mu$ on $\mathrm{Cvx}_{(o)}(\mathbb{R}^n)$ are immediate consequences of the definition of $\circledast$, it only remains to show that this map is continuous with respect to the topology induced by epi-convergence. Therefore, let $\varphi_k \in \mathrm{Cvx}_{(o)}(\mathbb{R}^n)$ be an epi-convergent sequence with limit $\varphi \in \mathrm{Cvx}_{(o)}(\mathbb{R}^n)$. In order to prove that
\begin{equation} \label{epiconvofconvol}
\varphi_k \circledast \mu \stackrel{\mathrm{epi}}{\rightarrow} \varphi \circledast \mu
\end{equation}
we proceed in three steps. First we claim that the Moreau envelope $e_t\varphi$ of $\varphi$ satisfies
\begin{equation} \label{moreauconvol}
e_t \varphi \circledast \mu \stackrel{\mathrm{epi}}{\rightarrow} \varphi \circledast \mu
\end{equation}
as $t \searrow 0$. Indeed, by (\ref{defcircledast}), Lemma \ref{propMoreau} (ii) and the monotone convergence theorem,
\[\lim_{t \searrow 0} (e_t \varphi \circledast \mu)(x) = \lim_{t \searrow 0}\! \int_{\mathbb{S}^{n-1}}\!\!\! e_t \varphi(\|x\| \vartheta_{x} v) \, d\mu(v)
= \!\int_{\mathbb{S}^{n-1}}\!\!\! \varphi(\|x\| \vartheta_{x} v) \, d\mu(v) = (\varphi \circledast \mu)(x)\]
for every $x \in \mathbb{R}^n \backslash\{o\}$, which, by Lemma \ref{lemepiconvequiv} (i), implies (\ref{moreauconvol}). In a second step, letting $k \rightarrow \infty$, we claim that for every $t > 0$,
\begin{equation} \label{moreauconvol1}
e_t \varphi_k \circledast \mu \stackrel{\mathrm{epi}}{\rightarrow} e_t\varphi \circledast \mu.
\end{equation}
To see this, we first note that $e_t \varphi_k$ is epi-convergent to $e_t \varphi$ as $k \rightarrow \infty$, by Proposition~\ref{propLegendre}~(iii) and
the fact that the Legendre transform maps infimal convolution to pointwise addition. Moreover, by Lemma \ref{propMoreau} (i),  $e_t \varphi_k$ and $e_t \varphi$ are both
finite and, hence, their epi-convergence is equivalent to uniform convergence on compact subsets of $\mathbb{R}^n$, by Lemma \ref{lemepiconvequiv} (ii), which in turn is preserved under the convolution $\circledast$ with the measure $\mu$.

Finally, we show that, letting $k \rightarrow \infty$,
\begin{equation} \label{whatwewant123}
(\varphi_k \circledast \mu)(x) \to (\varphi \circledast \mu)(x) \mbox{ for every } \left \{\!\! \begin{array}{l}  x \in \mathrm{int} \, \mathrm{dom}\, (\varphi \circledast \mu)\backslash\{o\}, \\
x \notin \mathrm{cl} \, \mathrm{dom}\, (\varphi \circledast \mu), \end{array} \right .
\end{equation}
which, by Lemma \ref{lemepiconvequiv} (i) and the fact that $\mathbb{R}^n$ without the boundary of $\mathrm{dom}\, (\varphi \circledast \mu)$ and the origin is a dense subset, concludes
the proof of (\ref{epiconvofconvol}) and the proposition.

To this end, first suppose that $x \not \in \mathrm{cl} \, \mathrm{dom}\, (\varphi \circledast \mu)$. Then there exists a closed ball $B$ such that $x \in B$ and
$B \cap \mathrm{cl} \, \mathrm{dom}\,(\varphi \circledast \mu) = \emptyset$. By (\ref{moreauconvol}) and Lemma \ref{lemepiconvequiv} (ii),
$e_t \varphi \circledast \mu$ converges to $\varphi \circledast \mu$ uniformly on $B$ as $t \searrow 0$. Thus, for every $c > 0$, we can find $t_0 > 0$
such that for every $t\leq t_0$, we have $e_t \varphi \circledast \mu \geq c$ on $B$. Similarly, by (\ref{moreauconvol1}) and Lemma \ref{lemepiconvequiv} (ii),
we can find $k_0 \in \mathbb{N}$ such that for all $k \geq k_0$, $e_t \varphi_k \circledast \mu > \frac{c}{2}$ on $B$.
Since $e_t \varphi_k \leq \varphi_k$, by Lemma \ref{propMoreau} (ii), and the convolution $\circledast$ is obviously monotone,
$\varphi_k \circledast \mu > \frac{c}{2}$ on $B$ for all $k \geq k_0$. Since $c > 0$ was arbitrary, we conclude that
$\varphi_k \circledast \mu$ converges to infinity uniformly on $B$, which proves (\ref{whatwewant123}) for $x \not \in \mathrm{cl} \, \mathrm{dom}\, (\varphi \circledast \mu)$.

Suppose now that $x \in \mathrm{int}\, \mathrm{dom}\, (\varphi \circledast \mu)$ is non-zero and fix a rotation $\vartheta_x \in \mathrm{SO}(n)$ such that $\vartheta_x \bar{e} = \frac{x}{\|x\|}$. Then there exists an open $\varepsilon$-ball $B_{\varepsilon}(x) \subseteq \mathrm{int}\, \mathrm{dom}\, (\varphi \circledast \mu)$ centered at $x$ . By (\ref{defcircledast}), $(\varphi_k \circledast \mu)(x)$ and $(\varphi \circledast \mu)(x)$ are determined by integrating the values of $\varphi_k$ and $\varphi$, respectively, over $\|x\| \mathbb{S}^{n-1}$. Therefore, we consider the compact set
\[C = (\|x\|\mathbb{S}^{n-1}) \cap (\mathbb{R}^n\backslash \mathrm{int}\,\mathrm{dom}\,\varphi).  \]
If $C$ is empty, then, by Lemma \ref{lemepiconvequiv} (ii), $\varphi_k$ converges uniformly to $\varphi$ on $\|x\|\mathbb{S}^{n-1}$ and, consequently, by (\ref{defcircledast}), we conclude that
$(\varphi_k \circledast \mu)(x) \to (\varphi \circledast \mu)(x)$.

Thus, assume that there exists some $y \in C$. Since $o \in \mathrm{int} \, \mathrm{dom}\, \varphi$,
there exists an open $\delta$-ball $B_{\delta}(o) \subseteq \mathrm{dom}\, \varphi$ centered at $o$. Since $\mathrm{dom}\,\varphi$ is convex,
each ray through $y \in C$ emanating from a point in $B_{\delta}(o)$ intersects the boundary of $\mathrm{dom}\,\varphi$ in exactly one point.
In particular, the parts of these rays starting at $y$ are completely contained in $\mathbb{R}^n\backslash \mathrm{dom}\,\varphi$. Hence, for every $y \in C$, there exists an open cone
$C_y$ with apex $y$ contained in $\mathrm{int}\,(\mathbb{R}^n\backslash \mathrm{dom}\,\varphi)$ and intersecting $R\,\mathbb{S}^{n-1}$, for any $R > \|x\|$, in an open cap whose diameter depends
only on $R$ and $\delta$.

Choosing $R = \|x\| + \frac{\varepsilon}{2}$, we have $x_0 = \frac{R}{\|x\|}x \in B_{\varepsilon}(x)$, that is, $\|x_0\| = R$ and $(\varphi \circledast \mu)(x_0) < \infty$ which imply, by (\ref{defcircledast}), that
\[\mu(\{u \in \mathbb{S}^{n-1}: \varphi(R\vartheta_xu)= \infty\}) = 0.  \]
Consequently, since $\varphi$ is infinite on $C_y \cap R\,\mathbb{S}^{n-1}$, we have
\[\mu(\{u \in \mathbb{S}^{n-1}: R\vartheta_xu \in C_y\}) = 0.  \]
Using that $R\vartheta_xu \in C_y$ if and only if $\|x\|\vartheta_xu \in \frac{\|x\|}{R}(C_y \cap R\,\mathbb{S}^{n-1}) \subseteq \|x\|\mathbb{S}^{n-1}$,
we infer that for every $y \in C$, there exists an open subset $U_y$ of $\mathbb{S}^{n-1}$ of $\mu$-measure zero such that $\|x\|\vartheta_xU_y$ is an open
neighborhood of $y$. The family $(\|x\| \vartheta_x U_y)_{y \in C}$ is an open cover of the compact set $C$, hence, there exists a finite subcover $(\|x\| \vartheta_x U_{y_i})_{i=1}^m$ and, by the sub-additivity of $\mu$, we have $\mu\left( U \right) = 0$ for $U = \bigcup_{i = 1}^m U_{y_i}$.

Since the compact set $C' = \|x\| \mathbb{S}^{n-1} \backslash \|x\| \vartheta_x U$ is disjoint from $\mathrm{bd}\,\mathrm{dom}\,\varphi$, $\varphi_k$ converges uniformly to $\varphi$ on $C'$, by Lemma \ref{lemepiconvequiv} (ii). For every $\widehat\varepsilon > 0$, we thus find $k_0$ such that for all $k \geq k_0$ we have $|\varphi_k(z) - \varphi(z)| \leq \widehat \varepsilon$ for all $z \in C'$. Consequently,
\begin{align*}
(\varphi_k \circledast \mu)(x) &= \int_{\mathbb{S}^{n-1}}\!\!\! \varphi_k(\|x\| \vartheta_x v) \, d\mu(v) = \int_{\mathbb{S}^{n-1} \backslash U}\!\!\! \varphi_k(\|x\| \vartheta_x v) \, d\mu(v)\\
                               &\leq \int_{\mathbb{S}^{n-1} \backslash U}\!\!\! \left (\varphi(\|x\| \vartheta_x v) + \widehat \varepsilon\,\right) d\mu(u) =
                               (\varphi \circledast \mu)(x) + \widehat \varepsilon\, \mu(\mathbb{S}^{n-1}),
\end{align*}
that is, $(\varphi_k \circledast \mu)(x)$ is finite for all $k \geq k_0$. Hence, we can infer that
\begin{align*}
|(\varphi_k \circledast \mu)(x) - (\varphi \circledast \mu)(x)| \leq \int_{\mathbb{S}^{n-1} \backslash U}\!\!\! |\varphi_k(\|x\| \vartheta_x v)-\varphi(\|x\| \vartheta_x v)| \, d\mu(v)
\leq \widehat \varepsilon\, \mu(\mathbb{S}^{n-1}),
\end{align*}
which completes the proof of (\ref{whatwewant123}).  \hfill $\blacksquare$

\vspace{0.4cm}

We are now in a position to complete the proof of Theorem \ref{def:phimuFirst}.

\vspace{0.3cm}

\noindent {\it Proof of Theorem \ref{def:phimuFirst}.} It follows from  (\ref{lccundcvx0}) and Proposition \ref{propforthm3} that $\Psi_{\mu}f$ is well defined for every $f \in \mathrm{LC_c}(\mathbb{R}^n)$. By (\ref{suppfctaspladd}) and (\ref{defcircledast}), $\Psi_{\mu}: \mathrm{LC_c}(\mathbb{R}^n) \to \mathrm{LC_c}(\mathbb{R}^n)$ is Asplund additive. By Lemma~\ref{asplproplem}~(iii), $f \star \mathbbm{1}_{\{y\}}$ coincides with the translate of $f \in \mathrm{LC_c}(\mathbb{R}^n)$ by $y \in \mathbb{R}^n$. Hence, using $\Psi_{\mu} \mathbbm{1}_{\{y\}} = \mathbbm{1}_{\Phi_{\mu}\{y\}}=\mathbbm{1}_{\{o\}}$, by Lemma~\ref{meanwidthminkendo} (iii), where $\Phi_{\mu}$ denotes the Minkowski endomorphism generated by $\mu$, we deduce that
\[\Psi_{\mu}(f \star \mathbbm{1}_{\{y\}})=\Psi_{\mu}f \star \mathbbm{1}_{\Phi_{\mu}\{y\}} = \Psi_{\mu}f \star \mathbbm{1}_{\{o\}} = \Psi_{\mu}f. \]
The commutativity of $\Psi_{\mu}$ with the action of $\mathrm{SO}(n)$ follows from the definitions of $\circledast$ and the support function $h(f,\cdot)$, and the fact that the Legendre transform commutes with the action of $\mathrm{SO}(n)$ on $\mathrm{Cvx}(\mathbb{R}^n)$. Continuity of $\Psi_{\mu}: \mathrm{LC_c}(\mathbb{R}^n) \to \mathrm{LC_c}(\mathbb{R}^n)$ is a consequence of Proposition \ref{propLegendre} (iii) and Proposition \ref{propforthm3}. The monotonicity of $\Psi_{\mu}$ follows from the monotonicity of support functions and that of $\circledast$.

By Proposition \ref{propLegendre} (iv) and the remark following (\ref{defcircledast}), $\Psi_{\mu} \mathbbm{1}_{K} = \mathbbm{1}_{\Phi_{\mu}K}$ for every $K \in \mathcal{K}^n$. Thus, $\mu$ is uniquely determined by $\Psi_{\mu}$ by Theorem \ref{weakmonMinkendoclass}. \hfill $\blacksquare$

\vspace{0.4cm}

In view of Theorem~\ref{weakmonMinkendoclass}, it is a natural question to ask whether Theorem \ref{def:phimuFirst} also holds for measures $\mu$ which are non-negative up to addition of a linear measure. In general, this is not possible as the integral of a (convex) function attaining the value $+\infty$ with respect to a signed measure is not well-defined. Restricting to the subspace of finite convex functions would overcome this issue, but it remains to prove that convoluting with $\mu$ preserves convexity.

\pagebreak 

The classification of additive (in a set-theoretic sense) maps on convex and log-concave functions has recently become the focus of intensive investigations \linebreak (see, e.g., \textbf{\cite{Alesker2019, Cavallina2015, Colesanti2017a, Colesanti2019, Colesanti2017b, Colesanti2019b, knoerr2020a, knoerr2020b}}). It is certainly an interesting open problem whether there exist monotone Asplund endomorphisms different from the ones provided by Theorem \ref{def:phimuFirst} and, if so, what additional properties characterize the endomorphisms from Theorem \ref{def:phimuFirst}.

We conclude this section with a functional analogue of Lemma \ref{meanwidthminkendo} (i).

\begin{lem} \label{meanwidthasplundendo} If $\mu$ is a non-negative zonal measure on $\mathbb{S}^{n-1}$ with center of mass at the origin, then
\begin{equation} \label{desired1742}
\int_{\mathbb{R}^n} h(\Psi_{\mu} f,x)\,d\gamma_n(x) =  \mu(\mathbb{S}^{n-1}) \int_{\mathbb{R}^n} h(f,x)\,d\gamma_n(x)
\end{equation}
for every $f \in \mathrm{LC_c}(\mathbb{R}^n)$.
\end{lem}

\noindent {\it Proof.} Using polar coordinates and the density $\psi_n$ of the Gaussian measure, yields
\begin{align*}
\int_{\mathbb{R}^n} h(\Psi_{\mu} f,x)\,d\gamma_n(x) &= n|B^n|\int_{0}^\infty \int_{\mathbb{S}^{n-1}}\!\!\! (h(f,\cdot) \circledast \mu)(ru)\, \psi_n(ru)\,r^{n-1} \,d \sigma(u)\, dr.
\end{align*}
By the $\mathrm{SO}(n)$ invariance of $\psi_n$, relation (\ref{desired1742}) follows, if we can show that
\begin{align} \label{chewbacca}
\int_{\mathbb{S}^{n-1}}\!\!\! (h \circledast \mu)(ru) \,d\sigma(u) = \mu(\mathbb{S}^{n-1})\int_{\mathbb{S}^{n-1}}\!\!\! h(ru) \,d\sigma(u)
\end{align}
for every $h \in \mathrm{Cvx}_{(o)}(\mathbb{R}^n)$ and every $r > 0$. To this end, first assume that $\mathrm{dom}\, h = \mathbb{R}^n$ and, hence, that $h$ is continuous.
Since $(h \circledast \mu)(ru) = (h(r\,\cdot\,) \ast \mu)(u)$ for every $u \in \mathbb{S}^{n-1}$, we obtain, as in the proof of Lemma \ref{meanwidthminkendo} (i), from the fact that spherical convolution is selfadjoint and Abelian for zonal measures, that
\begin{align*}
\int_{\mathbb{S}^{n-1}}\!\!\! (h \circledast \mu)(ru) \,d\sigma(u) =  \int_{\mathbb{S}^{n-1}} (h(r\,\cdot\,) \ast \sigma)(u)\,d\mu(u).
\end{align*}
By the $\mathrm{SO}(n)$ invariance of $\sigma$, $(h(r\,\cdot\,) \ast \sigma)(u)$ is independent of $u$ and (\ref{chewbacca}) follows.

For general $h \in \mathrm{Cvx}_{(o)}(\mathbb{R}^n)$, we use that the Moreau envelope $e_th$ of $h$ is convex and finite and converges monotonously to $h$,
by Lemma \ref{propMoreau}. By what we have shown above, (\ref{chewbacca}) holds for $e_th$ for every $t > 0$. Hence, by monotone convergence, we
conclude that (\ref{chewbacca}) holds generally. \hfill $\blacksquare$

\vspace{1cm}

\centerline{\large{\bf{ \setcounter{abschnitt}{4}
\arabic{abschnitt}. Proof of the main results}}}

\reseteqn \alpheqn \setcounter{theorem}{0}

\vspace{0.6cm}

In this section we first prove Theorem \ref{thm:ineqUrysohnMonotone} and deduce Theorem \ref{thm:ineqVolProdMonotone} from it. We then show that a further extension of Theorem \ref{thm:ineqVolProdMonotone} to all merely weakly monotone Minkowski endomorphisms is not possible. In order to prove Theorem \ref{thm:ineqVolProdFunc}, we will establish first a counterpart of Theorem \ref{thm:ineqUrysohnMonotone} for log-concave functions. We conclude this section by showing that each of the inequalities from Theorem \ref{thm:ineqVolProdFunc} is strictly stronger than the functional analogue of Urysohn's inequality.

\pagebreak

\noindent {\it Proof of Theorem \ref{thm:ineqUrysohnMonotone}.} Let $K \in \mathcal{K}^n$ have non-empty interior and note that the normalization $\Phi B^n = B^n$ ensures that there is equality
in both inequalities of (\ref{thm2inequ}) if $K$ is a Euclidean ball. By Theorem \ref{weakmonMinkendoclass}, the generating measure $\mu$ of $\Phi$ is non-negative and, by our normalization, $1 = h(\Phi B^n,\cdot) = \mu(\mathbb{S}^{n-1})$. Moreover, by Lemma~\ref{meanwidthminkendo}~(ii), $\Phi K$ contains the origin in its interior.

In order to establish the left hand inequality of (\ref{thm2inequ}), we use the polar coordinate formula for volume (\ref{polarformvol}), Jensen's inequality, and (\ref{defmeanwidth}) to obtain
\[\left ( \frac{|\Phi^\circ K|}{|B^n|}  \right )^{-1/n} = \left (\int_{\mathbb{S}^{n-1}}\!\!\! h(\Phi K,u)^{-n}d\sigma(u) \right )^{-1/n} \leq \int_{\mathbb{S}^{n-1}}\!\!\!  h(\Phi K,u)\,d\sigma(u)=\frac{w(\Phi K)}{2}.  \]
An application of Lemma \ref{meanwidthminkendo} (i) now yields the desired inequality. By the equality conditions for Jensen's inequality, equality holds here, and thus in the left hand side of (\ref{thm2inequ}),
if and only if $h(\Phi K,\cdot)$ is constant, that is, if and only if $\Phi K$ is a ball.

For the proof of the right hand inequality of (\ref{thm2inequ}), we may assume, by the translation invariance of both sides, that $\mathbf{s}(K) = o$, which implies $h(K,\cdot) > 0$ and $K^{\mathbf{s}}=K^\circ$. First, we use the polar coordinate formula for volume (\ref{polarformvol}), (\ref{zonalconv}), and Jensen's inequality (note that $\mu$ is a probability measure) to obtain
\begin{align*}
|\Phi^\circ K| &= \frac{1}{n} \int_{\mathbb{S}^{n-1}}\!\!\! h(\Phi K, u)^{-n} du = \frac{1}{n} \int_{\mathbb{S}^{n-1}}\! \left( \int_{\mathbb{S}^{n-1}}\!\!\! h(K, \vartheta_u v)\, d\mu(v)\! \right)^{-n}\! du \\
                   &\leq \frac{1}{n} \int_{\mathbb{S}^{n-1}}\!  \int_{\mathbb{S}^{n-1}}\!\!\! h(K, \vartheta_u v)^{-n}\, d\mu(v)\,  du.
\end{align*}
Since $\Phi$ and the polar map commute with $\mathrm{SO}(n)$ transforms, we may replace $K$ here by a rotated copy $\theta K$ and integrate over $\mathrm{SO}(n)$ with respect to the Haar measure, to arrive at
\begin{align} \label{eq:jensenMonotoneProof}
|\Phi^\circ K| = \int_{\mathrm{SO}(n)}\!\!\! |\Phi^\circ (\theta K)|\, d\theta \leq \frac{1}{n} \int_{\mathrm{SO}(n)}\! \int_{\mathbb{S}^{n-1}} \! \int_{\mathbb{S}^{n-1}}\!\!\! h(K, \theta^{-1} \vartheta_u v)^{-n}\,d\mu(v)\, du\,d \theta,
\end{align}
where we also used (\ref{sontranslsuppfct}) in the last step. By Fubini's theorem and the invariance of the Haar measure on $\mathrm{SO}(n)$,
\begin{align*}
|\Phi^\circ K| \leq \frac{1}{n} \int_{\mathbb{S}^{n-1}}\!\int_{\mathbb{S}^{n-1}} \! \int_{\mathrm{SO}(n)}\!\!\! h(K, \theta^{-1} u)^{-n}\,d \theta\,d\mu(v)\, du.
\end{align*}
Using the fact that $\mu(\mathbb{S}^{n-1})=1$ and again Fubini's theorem, we finally obtain the desired inequality,
\begin{align*}
|\Phi^\circ K| \leq \frac{1}{n} \int_{\mathbb{S}^{n-1}} \! \int_{\mathrm{SO}(n)}\!\!\! h(K, \theta^{-1} u)^{-n}\,d \theta\, du = \int_{\mathrm{SO}(n)}\!\!\! |\theta K^{\circ}|\, d\theta = |K^\circ|.
\end{align*}

By the above arguments, equality holds in the right hand inequality of (\ref{thm2inequ}) if and only if we have equality in (\ref{eq:jensenMonotoneProof}). By the equality condition of Jensen's inequality this is the case if and only if for almost every $u \in \mathbb{S}^{n-1}$ and almost every $\theta \in \mathrm{SO}(n)$ there exist $c_{u,\theta} \in \mathbb{R}^+$ such that
$h(K,\theta^{-1}\vartheta_uv) = c_{u,\theta}$ for $\mu$-a.e.\ $v \in \mathbb{S}^{n-1}$. Clearly, this is the case if and only if for every $\eta \in \mathrm{SO}(n)$ (by the continuity of $h(K,\cdot)$), there exist $c_{\eta} \in \mathbb{R}^+$ such that
\begin{equation} \label{geomsettmainprf}
h(K,\eta v) = c_{\eta} \quad \mbox{ for  $\mu$-a.e.\ $v \in \mathbb{S}^{n-1}$.}
\end{equation}
Let $A_{\eta} = \{v \in \mathbb{S}^{n-1}: h(K,\eta v) = c_{\eta}\} \subseteq \mathbb{S}^{n-1}$ and observe that $\mu(\mathbb{S}^{n-1}\backslash A_{\eta})=0$. Note that $A_{\eta}$ is closed, by the continuity of $h(K,\cdot)$, and therefore contains the support of $\mu$. If $\mu$ is discrete, by $\mathrm{SO}(n-1)$ invariance, it must coincide with the measure $\nu$ given by (\ref{discrzonal}) and, hence, $\Phi = \Delta$. Since $\mathrm{supp}\,\nu = \{-\bar{e},\bar{e}\}$, it follows immediately from (\ref{geomsettmainprf}) that $h(K,\cdot)$ takes the same value on antipodal points, that is, $K$ is origin-symmetric.

It remains to be shown that if $\mu$ is not discrete, then (\ref{geomsettmainprf}) holds if and only if $h(K,\cdot)$ is constant on $\mathbb{S}^{n-1}$, or equivalently, if $K$ is an origin-symmetric Euclidean ball. Since $\mu$ is non-zero and not discrete, there exists  $w \in \mathrm{supp}\, \mu \backslash \{-\bar{e},\bar{e}\}$. By the $\mathrm{SO}(n-1)$-invariance of $\mu$, the entire parallel subsphere
orthogonal to $\bar e$ through $w$ is contained in $\mathrm{supp}\, \mu$. Hence, by (\ref{geomsettmainprf}), $h(K,\eta v) = c_{\eta}$ for every $\eta \in \mathrm{SO}(n)$ and all $v$ in this subsphere.
Choosing $\eta'$ such that this subsphere and a copy of it rotated by $\eta'$ intersect, we see that the value of $h(K,\cdot)$ at the intersection is given by
$c_{\mathrm{id}}$ and $c_{\eta'}$, thus, these values must be equal. By repeating this argument finitely many times, we can reach every point on $\mathbb{S}^{n-1}$
implying that $h(K,\cdot)$ is constant as desired.  \hfill $\blacksquare$

\vspace{0.4cm}

Theorem \ref{thm:ineqVolProdMonotone} is now an easy consequence of the right-hand inequality of (\ref{thm2inequ}).

\vspace{0.3cm}

\noindent {\it Proof of Theorem \ref{thm:ineqVolProdMonotone}.}  Let $K \in \mathcal{K}^n$ have non-empty interior and assume w.l.o.g.\ that $\Phi B^n = B^n$.
Then, by the right-hand inequality of (\ref{thm2inequ}) and the Blaschke--Santal\'o inequality (\ref{BSinequ}), it follows that
\[|K||\Phi^{\circ}K| \leq |K||K^{\mathbf{s}}| \leq |B^n|^2 = |B^n||\Phi^{\circ}B^n|.  \]
The equality $|K||\Phi^{\circ}K|=|B^n|^2$ holds if and only if equality holds both in the right-hand inequality of (\ref{thm2inequ})
and the Blaschke--Santal\'o inequality (\ref{BSinequ}), that is, if and only if $\Phi = \Delta$ and $K$ is an ellipsoid or if $K$ is a Euclidean ball. \hfill $\blacksquare$

\vspace{0.4cm}

Next, we want to show that an extension of Theorem \ref{thm:ineqVolProdMonotone} to all merely weakly monotone Minkowski endomorphisms is impossible.

\begin{theorem} \label{Junbounded} For every $n \geq 2$, the volume product $|K||\mathrm{J}^{\circ}K|$ is unbounded for the weakly monotone Minkowski endomorphism
$\mathrm{J}: \mathcal{K}^n \rightarrow \mathcal{K}^n$, $\mathrm{J}K = K - s(K)$.
\end{theorem}

\noindent {\it Proof.} We begin with dimension $n = 2$, where for $c > 0$, we consider the triangle $K_c \in \mathcal{K}^2$ of unit volume defined by
$K_c = \text{conv}\left\{\left(c, 0\right), \left(0,\frac{1}{c}\right), \left(0, -\frac{1}{c}\right)\right\}$. Then for every $0 < t < c$, the polar body
$(K_c - (t,0))^\circ$ is again a triangle given by
\[(K_c-(t,0))^\circ = \text{conv}\left\{ \left (\frac{1}{c-t} , \frac{c^2}{c-t}\right), \left (\frac{1}{c-t}, -\frac{c^2}{c-t}\right), \left (-\frac{1}{t}, 0\right)   \right \}. \]
Thus, a short calculation yields the volume formula,
\begin{align} \label{leia1742}
\left | \left(K_c - \left(t, 0\right) \right)^\circ \right| = \frac{c^3}{t(c-t)^2}.
\end{align}
Due to the axial symmetry of $K_c$, its Steiner point $s(K_c)$ lies on the $x$-axis and its coordinates are easily calculated to
\begin{align*}
s(K_c) = \left( \frac{c}{\pi}\arctan{c^2}, 0 \right)\!.
\end{align*}
Plugging this into (\ref{leia1742}), we obtain
\begin{align*}
|K_c|\left| \left(K_c - s(K_c) \right)^\circ \right| = \frac{\pi}{\arctan(c^2)(1-\frac{1}{\pi}\arctan(c^2))^2},
\end{align*}
which tends to infinity as $c$ tends to zero (by a computation involving only elementary calculus).

For arbitrary $n \geq 3$, we consider the body of revolution $L_c \in \mathcal{K}^n$, obtained by rotating the body $K_c$ around the $e_1$-axis of $\mathbb{R}^n$.
The volume of $L_c$ can be easily calculated and is given by
\begin{align} \label{volLc}
|L_c| = \frac{|B^{n-1}|}{nc^{n-2}},
\end{align}
where $|B^{n-1}|$ denotes the $(n-1)$-dimensional Hausdorff measure of the Euclidean unit ball in $\mathbb{R}^{n-1}$.

Since for every $K \in \mathcal{K}^n$ containing the origin in its interior and any subspace $H \subseteq \mathbb{R}^n$, we have
$K^ \circ | H = (K \cap H)^\circ$, where the polar body on the right hand side is taken in the subspace $H$, it follows, by taking $H$ a $2$-dimensional
subspace containing $e_1$, that for every $0 < t < c$, $(L_c - te_1)^\circ$ is a body of revolution obtained by rotating the triangle $(K_c-te_1)^\circ$ around the
$e_1$ axis. Consequently, we obtain the volume formula,
\begin{align*}
|(L_c - t e_1)^\circ| = \frac{|B^{n-1}|}{n} \frac{c^{2n-1}}{t(c-t)^{n}}
\end{align*}
and, from this and (\ref{volLc}),
\begin{align*}
|L_c||(L_c-t e_1)^\circ| = \left(\frac{|B^{n-1}|}{n}\right)^2 \frac{c^{n+1}}{t(c-t)^{n}}.
\end{align*}
Letting $t = c g(c)$, where $g(c)$ depends only on $c$ and satisfies $0 < g(c) < 1$, this reduces to
\begin{align*}
|L_c||(L_c-cg(c)e_1)^\circ| = \left(\frac{|B^{n-1}|}{n}\right)^2 \frac{1}{g(c)(1-g(c))^{n}}
\end{align*}
which clearly tends to infinity if $g(c)$ tends to zero as $c$ tends to zero. It remains to be shown that the $e_1$ coordinate of the Steiner point of $L_c$ is
of the form $cg(c)$ such that $\lim_{c \to 0}g(c) = 0$ (note that by the rotational symmetry of $L_c$ all other coordinates of $s(L_c)$ are zero).

\pagebreak

Since $h(L_c,u) = h(K_c,(u \cdot e_1) e_1 + \sqrt{1-(u\cdot e_1)^2} e_2)$
for every $u \in \mathbb{S}^{n-1}$, we obtain from (\ref{defsteinerpkt}) by integration in cylindrical coordinates (cf.\ \textbf{\cite[\textnormal{Lemma~1.3.1}]{Groemer1996}}),
\begin{align*}
 s(L_c)\! \cdot e_1 = \frac{1}{|B^n|}\! \int_{\mathbb{S}^{n-1}}\!\!\!\!\!\!\! h(L_c,u)(u \cdot e_1)du = n\!\! \int_{-1}^1\!\!\! h(K_c,\zeta e_1 + \sqrt{1-\zeta^2} e_2) \zeta (1-\zeta^2)^{\frac{n-3}{2}}d\zeta.
\end{align*}

Inserting the explicit expression for the support function of $K_c$, we see that
\begin{align*}
s(L_c) \cdot e_1 & = n \left( \int_{-1}^{\frac{1}{\sqrt{1+c^4}}} \frac{\zeta}{c}(1-\zeta^2)^{\frac{n-1}{2}}\,d\zeta + \int_{\frac{1}{\sqrt{1+c^4}}}^1 c \zeta^2 (1-\zeta^2)^{\frac{n-3}{2}}\,d\zeta\right) \\ & = c \underbrace{n\left(\frac{-c^{2n}}{(n+1)(1+c^4)^\frac{n+1}{2}} + \int_{\frac{1}{\sqrt{1+c^4}}}^1 \zeta^2 (1-\zeta^2)^{\frac{n-3}{2}}\,d\zeta \right)}_{= g(c)}
\end{align*}
which is of the desired form. \hfill $\blacksquare$

\vspace{0.4cm}

Theorem \ref{Junbounded} raises the interesting problem whether there exist weakly monotone or even non-monotonic Minkowski endomorphisms $\Phi$ (different from multiples of $\mathrm{J}$ and $-\mathrm{J}$) such that their volume product $|K||\Phi^\circ K|$ is unbounded. In view of Theorem \ref{thm:ineqVolProdMonotone}, it is also natural to ask whether there exist Minkowski endomorphisms wich are maximized by convex bodies of non-empty
interior which may be different from Euclidean balls. Partial answers to this question were very recently obtained in \textbf{\cite{hofstaetteretal2021}} by methods different to the ones used in this article.

Let us also comment on our general assumption that $n \geq 3$. An analogue of Theorem~\ref{weakmonMinkendoclass} for $n=2$ was already obtained by Schneider~\textbf{\cite{Schneider1974}} in 1974 showing also that all Minkowski endomorphisms in $\mathbb{R}^2$ are weakly monotone. The key difference in this case is the commutativity of $\mathrm{SO}(2)$, which implies that for every monotone Minkowski endomorphism $\Phi$, $\Phi K$ is the limit of combinations $\lambda_1 \vartheta_1 K + \dots + \lambda_m \vartheta_m K$, where $\vartheta_i \in \mathrm{SO}(2)$ and $\lambda_i > 0$ such that $\sum_{i=1}^m \lambda_i$ is fixed. The well-known inequality
$ \left| \left( s_1 K_1 + s_2 K_2\right)^\circ \right| \leq s_1 |K_1^\circ| + s_2 |K_2^\circ|$ for $K_1, K_2 \in \mathcal{K}^2$ with non-empty interior and $s_1, s_2 > 0$ such $s_1 + s_2 = 1$, directly implies an analogue of Theorem~\ref{thm:ineqUrysohnMonotone}.

\vspace{0.2cm}

We turn now to inequalities for log-concave functions and begin with a functional analogue of Theorem \ref{thm:ineqUrysohnMonotone} from which we will subsequently
deduce Theorem \ref{thm:ineqVolProdFunc}.

\begin{theorem} \label{functanalogthm2} Let $\mu$ be an $\mathrm{SO}(n-1)$ invariant probability measure on $\mathbb{S}^{n-1}$ with center of mass at the origin. If $f\in \mathrm{LC_c}(\mathbb{R}^n)$ such that $\int_{\mathbb{R}^n} f\,dx > 0$, then
\begin{align} \label{maininequchainfct}
\int_{\mathbb{R}^n}(\Psi_{\sigma} f)^\circ(x)\, dx \leq \int_{\mathbb{R}^n}(\Psi_{\mu} f)^\circ(x)\, dx \leq   \int_{\mathbb{R}^n} f^\circ(x)\, dx.
\end{align}
There is equality in the left hand inequality if and only if $\Psi_{\mu}f$ is radially symmetric. Equality in the right hand inequality holds if and only if $f$ is even and $\Psi_{\mu} = \Delta_{\star}$ or if $f$ is radially symmetric.
\end{theorem}

\noindent {\it Proof.} First note that for $f \in \mathrm{LC_c}(\mathbb{R}^n)$, we always have $\int_{\mathbb{R}^n} f\,dx < \infty$, by (\ref{coercconecond}). In order to establish the left hand inequality of (\ref{maininequchainfct}), we use (\ref{defpolarf}), polar coordinates and Jensen's inequality to obtain
\begin{align}
\int_{\mathbb{R}^n}\!(\Psi_{\mu} f)^\circ(x)\, dx & = n|B^n| \int_0^{\infty}\!\! \int_{\mathbb{S}^{n-1}}\!\!\! \exp\!\left(-(h(f,\cdot) \circledast \mu)(ru)\right )\,r^{n-1}\,d\sigma(u)\,dr \nonumber \\
& \geq n|B^n| \int_0^{\infty}\!\!  \exp\!\left (\!-\!\int_{\mathbb{S}^{n-1}}\!\!\! (h(f,\cdot) \circledast \mu)(ru)\,d\sigma(u)\!\right)r^{n-1}\,dr. \label{test17}
\end{align}
To be precise, for the application of Jensen's inequality, we require the function $(h(f,\cdot) \circledast \mu)(r\,\cdot\,)$ to be $\sigma$-integrable. However, if this is not the case, then its integral is $+\infty$ and inequality (\ref{test17}) still holds.

From an application of (\ref{chewbacca}) to the inner integral in (\ref{test17}) and the $\mathrm{SO}(n)$ invariance of $\sigma$, we conclude that
\begin{align*}
\int_{\mathbb{R}^n}\!(\Psi_{\mu} f)^\circ(x)\, dx & \geq n|B^n| \int_0^{\infty}\!\!  \exp\!\left (\!-\!\int_{\mathbb{S}^{n-1}}\!\!\! h(f,ru) \,d\sigma(u)\!\right)r^{n-1}\,dr \\
& = n|B^n| \int_0^{\infty}\!\!  \exp\!\left (-(h(f,\cdot) \circledast \sigma)(rv) \right ) r^{n-1}\,dr
\end{align*}
for an arbitrary $v \in \mathbb{S}^{n-1}$. Finally, using that $(h(f,\cdot) \circledast \sigma)(rv)$ does not depend on $v$, we arrive at the left hand inequality of (\ref{maininequchainfct}),
\begin{align*}
\int_{\mathbb{R}^n}\!(\Psi_{\mu} f)^\circ(x)\, dx & \geq n|B^n|\! \int_0^{\infty}\!\!\int_{\mathbb{S}^{n-1}}\!\!\!
\exp\!\left (-(h(f,\cdot) \circledast \sigma)(rv) \right ) r^{n-1}\,d\sigma(v)\,dr \\ & = \int_{\mathbb{R}^n}\!(\Psi_{\sigma} f)^\circ(x) \,dx.
\end{align*}
If equality holds in the left hand inequality of (\ref{maininequchainfct}), then we must have equality in (\ref{test17}) which implies by the equality condition of Jensen's inequality (including the case of non-$\sigma$-integrability) that for almost every $r > 0$ there exists $c_r \in (-\infty,\infty]$ such that
\begin{equation} \label{diffprf17}
h(\Psi_{\mu}f,rv) = (h(f,\cdot) \circledast \mu)(rv) = c_r \quad \mbox{ for $\sigma$-a.e.\ $v \in \mathbb{S}^{n-1}$.}
\end{equation}
Note that, by continuity, (\ref{diffprf17}) yields that $h(\Psi_{\mu}f,\cdot)$ is constant on every sphere contained in $\mathrm{int}\,\mathrm{dom}\,h(\Psi_\mu f,\cdot)$. Next, we want to show that this domain is a ball.

If $c_r < \infty$ for some $r > 0$, then the lower semi-continuity of $h(\Psi_{\mu}f,\cdot)$ implies that $h(\Psi_{\mu}f,rv)$ is finite for \emph{every} $v \in \mathbb{S}^{n-1}$.
In particular, $r\mathbb{S}^{n-1} \subseteq \mathrm{dom}\,h(\Psi_{\mu}f,\cdot)$ which by the convexity of this domain yields $rB^n \subseteq \mathrm{dom}\,h(\Psi_{\mu}f,\cdot)$.
This implies that (\ref{diffprf17}) holds for every $r' \leq r$, by continuity, and $c_{r'} < \infty$. Thus, the set of all $r > 0$ such that (\ref{diffprf17}) holds with $c_r < \infty$ is an interval, that is, there
exists $R > 0$ such that
\[c_r \left \{\!\!\begin{array}{ll} < \infty & \mbox{ for all } r < R, \\ = \infty & \mbox{ for all } r > R \mbox{ for which (\ref{diffprf17}) holds.} \end{array} \right .   \]
In order to conclude that $\mathrm{int}\,\mathrm{dom}\,h(\Psi_{\mu}f,\cdot)$ is a ball, it remains to show that for every $r > R$, $h(\Psi_{\mu}f,\cdot)$ is infinite on $r\mathbb{S}^{n-1}$.
To this end, let $x \in r\mathbb{S}^{n-1}$ and assume that $h(\Psi_{\mu}f,x) < \infty$. Since $\mathrm{dom}\,h(\Psi_{\mu}f,\cdot)$ is convex and contains an open ball centered at the
origin, the convex hull of $x$ and this ball is contained in $\mathrm{dom}\,h(\Psi_{\mu}f,\cdot)$. However, this convex hull must contain an open neighborhood of
$r'\,\mathbb{S}^{n-1}$ for some $r > r' > R$ for which (\ref{diffprf17}) holds, which contradicts $c_{r'} = \infty$.

Finally, since $\mathrm{int}\,\mathrm{dom}\,h(\Psi_\mu f,\cdot)$ is a ball, by the comment following (\ref{diffprf17}), $h(\Psi_\mu f,\cdot)$ is radially symmetric
on this ball. As a convex function which is radially symmetric (and, thus, depends only on one variable) on the interior of its domain must be radially symmetric everywhere, we conclude that $h(\Psi_\mu f,\cdot)$ is radially symmetric on $\mathbb{R}^n$. This concludes the proof of the equality conditions for the left hand inequality of (\ref{maininequchainfct}).

For the proof of the right hand inequality of (\ref{maininequchainfct}), we use (\ref{defpolarf}), (\ref{defcircledast}), and Jensen's inequality (note that $\mu$ is a probability measure) to obtain
\begin{align} \label{jarjarbinks}
\int_{\mathbb{R}^n}\!(\Psi_\mu f)^\circ(x)\,dx \leq \int_{\mathbb{R}^n} \int_{\mathbb{S}^{n-1}}\!\!\!
 \exp\left(-h(f,\|x\| \vartheta_x v)\right)\,d\mu(v)\,dx.
\end{align}
As in the first part of this proof, for the application of Jensen's inequality, we require $h(f,\|x\|\vartheta_x\,\cdot)$ to be $\mu$-integrable.
However, if this is not the case, then the left hand side of (\ref{jarjarbinks}) is zero and inequality (\ref{jarjarbinks}) still holds.

Since $\Psi_{\mu}$ and the polar map commute with $\mathrm{SO}(n)$ transforms, we may replace $f$ in (\ref{jarjarbinks})
 by a rotated copy $\theta f$ and integrate over $\mathrm{SO}(n)$ with respect to the Haar measure, to arrive at
\begin{align} \label{thisistheinequ}
\int_{\mathbb{R}^n}\!(\Psi_\mu f)^\circ(x)\,dx \leq \int_{\mathrm{SO}(n)}\! \int_{\mathbb{R}^n} \int_{\mathbb{S}^{n-1}}\!\!\!
 \exp\left(-h(f,\|x\| \theta^{-1}\vartheta_x v)\right)\,d\mu(v)\,dx\,d\theta.
\end{align}
Since we integrate non-negative functions, we may apply Fubini's theorem twice, the invariance of the Haar measure on $\mathrm{SO}(n)$,
and the fact that $\mu(\mathbb{S}^{n-1})=1$, to obtain the desired inequality
\begin{align*}
\int_{\mathbb{R}^n}\!(\Psi_\mu f)^\circ(x)\,dx & \leq \int_{\mathbb{R}^n} \int_{\mathbb{S}^{n-1}}\int_{\mathrm{SO}(n)}\!\!\!
 \exp\left(-h(f,\|x\| \theta^{-1}\vartheta_x v)\right)d\theta\,d\mu(v)\,dx \\
  &=  \int_{\mathbb{R}^n}\int_{\mathrm{SO}(n)}\!\!\! \exp\left(-h(f, \theta^{-1} x)\right)\,d\theta\,dx \\
  &=  \int_{\mathrm{SO}(n)}\! \int_{\mathbb{R}^n}\! f^\circ(\theta^{-1}x)\,dx\,d\theta = \int_{\mathbb{R}^n} f^\circ(x)\,dx.
 \end{align*}

By the above arguments, equality holds in the right hand inequality of (\ref{maininequchainfct}) if and only if we have equality in (\ref{thisistheinequ})
which implies by the equality condition of Jensen's inequality (including the case of non-$\mu$-integrability) that for
almost every $\theta \in \mathrm{SO}(n)$ and almost every $x \in \mathbb{R}^n$ there exist constants $c_{\theta,x} \in (-\infty,\infty]$ such that
\begin{equation} \label{thisisthecond}
h(f,\|x\|\theta^{-1}\vartheta_x v) = c_{\theta,x} \quad \mbox{ for $\mu$-a.e.\ } v \in \mathbb{S}^{n-1}.
\end{equation}

As in the proof of Theorem \ref{thm:ineqUrysohnMonotone}, if $\mu$ is discrete it must coincide with the measure $\nu$ given by (\ref{discrzonal}).
Thus, $\Psi_{\mu} = \Delta_{\star}$ and, since $\mathrm{supp}\,\nu = \{-\bar{e},\bar{e}\}$, (\ref{thisisthecond}) reduces to the existence of constants
$c_{\theta,x} \in (-\infty,\infty]$ such that for almost every $\theta \in \mathrm{SO}(n)$ and almost every $x \in \mathbb{R}^n$,
 \begin{align*}
h(f,\theta^{-1}x) = c_{\theta, x} = h(f,-\theta^{-1}x).
 \end{align*}
Consequently, the interior of the domain of $h(f,\cdot)$ must be origin-symmetric and $h(f,\cdot)$ must be even on it (by continuity, $h(f,\cdot)$ must attain the same value
on \emph{all} antipodal points in $\mathrm{int}\,\mathrm{dom}\,h(f,\cdot)$). By now considering the restriction of $h(f,\cdot)$ to lines through the origin
and using the extendibility of convex, lower semi-continuous functions of one variable, we conclude that $h(f,\cdot)$ must be even on all of $\mathbb{R}^n$.

If $\mu$ is not discrete, we first want to show that (\ref{thisisthecond}) implies that $\mathrm{int}\, \mathrm{dom}\,h(f,\cdot)$ is an open ball centered at the
origin. To this end, note that it follows from (\ref{thisisthecond}) that for all $r$ from a dense subset of $(0,\infty)$ and almost every $\eta \in \mathrm{SO}(n)$, there exist constants
$c_{r,\eta} \in (-\infty,\infty]$ such that
\begin{align} \label{thisisthenewcond}
h(f,r \eta v) = c_{r,\eta} \quad \text{ for $\mu$-a.e.\ } v \in \mathbb{S}^{n-1}.
\end{align}
If $c_{r,\eta}< \infty$ for some $r > 0$ and $\eta \in \mathrm{SO}(n)$, then the lower semi-continuity of $h(f,\cdot)$ implies that
$h(f,r\eta v) \leq c_{r,\eta}< \infty$ for all $v \in \mathrm{supp}\,\mu$. If on the other hand $c_{r,\eta} = \infty$, then the lower semi-continuity
of $h(f,\cdot)$ implies that $v \mapsto h(f,r \eta v)$ cannot be bounded on any open subset of $\mathbb{S}^{n-1}$ intersecting $\mathrm{supp}\, \mu$.

For suitable $\delta > 0$, let $B_{\delta}$ denote an open origin-symmetric $\delta$-ball in $\mathrm{dom}\,h(f,\cdot)$ such that its closure is still contained
in $\mathrm{int}\,\mathrm{dom}\,h(f,\cdot)$. Next, choose an arbitrary $x \in \mathrm{int}\,\mathrm{dom}\,h(f,\cdot) \backslash \mathrm{cl}\, B_{\delta}$ and
let $r \in [\delta,\|x\|)$. Then the set $\mathrm{conv}\{x, \mathrm{cl}\,B_{\delta}\}$ is contained in $\mathrm{int}\,\mathrm{dom}\,h(f,\cdot)$ and, thus,
$h(f,\cdot)$ is bounded on it. Define the open spherical caps $C_x^r$ as $\mathrm{conv}\{x, B_{\delta}\} \cap r\mathbb{S}^{n-1}$.

In the following, let $d$ denote the geodesic distance on $\mathbb{S}^{n-1}$ and let $d_r$ denote the geodesic distance on $r\mathbb{S}^{n-1}$
normalized such that $d_r(u,v) = d(\frac{u}{r}, \frac{v}{r})$ for any $u,v\in r\mathbb{S}^{n-1}$. Since $\mu$ is not discrete and has center of mass at the origin, there exists $t_0 \in [0,1)$ such that
$H_{\bar e, t_0} \cap \mathbb{S}^{n-1} \subseteq \mathrm{supp}\, \mu$, where $H_{\bar e, t_0} = \{y \in \mathbb{R}^n: \bar{e} \cdot y = t_0\}$.
Let $\alpha$ be the maximal geodesic distance of two points in $H_{\bar e, t_0} \cap \mathbb{S}^{n-1}$.

Choose now $r_0$ from the dense subset of all $r \in [\delta,\|x\|)$ such that (\ref{thisisthenewcond}) holds for almost all $\eta \in \mathrm{SO}(n)$.
Then for $x_0 \in r_0\mathbb{S}^{n-1}$ with  $d_{r_0}(x_0,C_x^{r_0}) < \frac{\alpha}{2}$ and $\varepsilon > 0$, we consider the set
\begin{align*}
A_{x_0, \varepsilon} = \{\eta \in \mathrm{SO}(n): r_0 \eta v_1 \in C_x^{r_0}, d_{r_0}(r_0 \eta v_2, x_0) < \varepsilon \text{ for some } v_1, v_2 \in H_{\bar e, t_0} \cap \mathbb{S}^{n-1} \}.
\end{align*}
Clearly, $A_{x_0, \varepsilon}$ is open and non-empty, hence, there exists $\eta_0 \in A_{x_0, \varepsilon}$ such that (\ref{thisisthenewcond}) holds for
$\eta_0$ and $r_0$. Since $C_x^{r_0}$ is open and $h(f,\cdot)$ is bounded on $C_x^{r_0}$, we have $c_{r_0,\eta_0} < \infty$ (as we have seen above).
In particular, $h(f,r_0 \eta_0 v_2) < \infty$, that is, $r_0 \eta_0 v_2 \in \mathrm{dom}\,h(f,\cdot)$ for some $v_2 \in H_{\bar e, t_0} \cap \mathbb{S}^{n-1}$ such that
$d_{r_0}(r_0 \eta_0 v_2, x_0) < \varepsilon$ and there exists $v_1 \in H_{\bar e, t_0} \cap \mathbb{S}^{n-1}$ with $r_0 \eta_0 v_1 \in C_x^{r_0}$.

Since $x_0 \in r_0\mathbb{S}^{n-1}$ such that $d_{r_0}(x_0,C_x^{r_0}) < \frac{\alpha}{2}$ and $\varepsilon > 0$ were arbitrary, $h(f,\cdot)$ is finite on a dense subset
of $U_{\alpha/2}(C_x^{r_0})$, the set of all points on $r_0 \mathbb{S}^{n-1}$ whose distance $d_{r_0}$ to $C_x^{r_0}$ is less than $\frac{\alpha}{2}$.
Taking any $r' < r_0$, this implies that $U_{\alpha/2}(C_x^{r'}) \subseteq \mathrm{dom}\,h(f,\cdot)$. Indeed, for every $y \in U_{\alpha/2}(C_x^{r'})$
we may find points $x_1, \dots, x_n$ in $U_{\alpha/2}(C_x^{r_0}) \cap \mathrm{dom}\,h(f,\cdot)$ such that
$y \in \mathrm{conv}\{0, x_1, \dots, x_n\} \subseteq \mathrm{dom}\,h(f,\cdot)$.
Since $r_0$ can be chosen from a dense subset of $[\delta,\|x\|)$, we have shown that for every $x \in \mathrm{int}\,\mathrm{dom}\,h(f,\cdot)$ and \emph{every} $r \in [\delta,\|x\|)$,
the set $U_{\alpha/2}(C_x^r)$ is contained in $\mathrm{dom}\,h(f,\cdot)$. In particular, $U_{\alpha/2}(\frac{r}{\|x\|}x) \subseteq \mathrm{dom}\,h(f,\cdot)$.

Finally, if $\mathrm{int}\,\mathrm{dom}\,h(f,\cdot)$ is not a centered open ball, then there exists a sequence $x_k \in \mathrm{int}\,\mathrm{dom}\,h(f,\cdot)$ converging to $x \in \mathrm{bd}\,\mathrm{dom}\,h(f,\cdot)$ with $\|x_k\| > \|x\| > 0$ for all~$k$ (take, e.g., for $x$ any boundary point that is touched non-radially by a closed ball in $\mathrm{int}\,\mathrm{dom}\,h(f,\cdot)$).

Since $x_k \rightarrow x$, $\frac{\|x\|}{\|x_k\|}x_k$ converges to $x$ as well. Hence, there exists $k_0 \in \mathbb{N}$ such that
$d_{\|x\|}\!\left (\frac{\|x\|}{\|x_{k_0}\|}x_{k_0}, x \right ) < \frac{\alpha}{4}$, that is,
\[ x \in U_{\alpha/4}\!\left (\frac{\|x\|}{\|x_{k_0}\|}x_{k_0}\right) \subseteq \bigcup_{r \in (\delta, \|x_{k_0}\|)}U_{\alpha/2}\!\left (\frac{r}{\|x_{k_0}\|}x_{k_0}\right)
\subseteq \mathrm{int}\,\mathrm{dom}\,h(f,\cdot) \]
which is a contradiction.

Knowing that $\mathrm{int}\,\mathrm{dom}\,h(f,\cdot)$ is a centered open ball, (\ref{thisisthenewcond}) combined with the continuity of $h(f,\cdot)$ on the interior of its
domain, implies, as in the final paragraph of the proof of Theorem \ref{thm:ineqUrysohnMonotone}, that $h(f,\cdot)$
is radially symmetric on $\mathrm{int}\,\mathrm{dom}\,h(f,\cdot)$. Noting again that a convex function which is radially symmetric on the interior of its domain must be radially symmetric everywhere, we infer that $h(f,\cdot)$ is radially symmetric on all of $\mathbb{R}^n$. Since the Legendre transform commutes with
the action of $\mathrm{SO}(n)$, $f$ must be radially symmetric itself. As $\Psi_{\mu}f = f$ for any radially symmetric $f \in \mathrm{LC_c}(\mathbb{R}^n)$, by (\ref{defcircledast}), this
concludes the proof of the theorem. \hfill $\blacksquare$

\vspace{0.3cm}

The same way Theorem \ref{thm:ineqVolProdMonotone} was a simple consequence of Theorem \ref{thm:ineqUrysohnMonotone} and (\ref{BSinequ}), we can now deduce Theorem \ref{thm:ineqVolProdFunc} easily from Theorems~\ref{functanalogthm2} and \ref{generalfctBS}.

\vspace{0.3cm}

\noindent {\it Proof of Theorem \ref{thm:ineqVolProdFunc}.}  By the translation-invariance of $\Psi_{\mu}$, we have $\Psi_{\mu}f = \Psi_{\mu}\tilde{f}$, where, as in Theorem \ref{generalfctBS},
$\tilde{f}(x) = f(x-\mathrm{cent}\,f)$. Thus, by the right hand inequality of (\ref{maininequchainfct}) and Theorem \ref{generalfctBS}, it follows that
\begin{align*}
\int_{\mathbb{R}^n}\! f(x) \,dx \int_{\mathbb{R}^n}\!(\Psi_{\mu} f)^\circ(x)\,dx  &= \int_{\mathbb{R}^n}\! f(x)\,dx \int_{\mathbb{R}^n}\!(\Psi_{\mu} \tilde f)^\circ(x)\,dx \\
                                                        &\leq \int_{\mathbb{R}^n}\! f(x)\,dx \int_{\mathbb{R}^n}\!\tilde{f}^\circ(x)\,dx \leq (2\pi)^n.
\end{align*}
Equality holds in (\ref{functinequAsplendos}) if and only if equality holds both in the right hand inequality of (\ref{maininequchainfct}) and in Theorem \ref{generalfctBS}, that is, if and only if $\Psi_{\mu} = \Delta_{\star}$ and $f$ is a Gaussian or if $f$ is proportional to a translation of the standard Gaussian. \hfill $\blacksquare$

\vspace{0.3cm}

Let us remark at this point again that Theorem \ref{thm:ineqVolProdMonotone} can be recovered from Theorem~\ref{thm:ineqVolProdFunc} in an asymptotically optimal form. More precisely, choosing $f = \mathbbm{1}_K$ for $K \in \mathcal{K}^n$ with non-empty interior in Theorem \ref{thm:ineqVolProdFunc}, inequality (\ref{functinequAsplendos}) becomes,
\[(2\pi)^n \geq |K| \int_{\mathbb{R}^n} \mathbbm{1}_{\Phi_{\mu}K}^\circ(x)\,dx = |K| \int_{\mathbb{R}^n} \exp\!\left (-\|x\|_{\Phi_{\mu}^\circ K} \right )dx = n!|K||\Phi_{\mu}^\circ K|,  \]
where we have used that $\Psi_{\mu}\mathbbm{1}_K = \mathbbm{1}_{\Phi_{\mu}K}$, the definition of the polar map, and (\ref{volform17}) with $p = 1$. Since the assumption that $\mu$ is a probability measure is equivalent to the normalization $\Phi_{\mu}B^n = B^n$, we obtain
\[|K||\Phi_{\mu}^\circ K| \leq \frac{(2\pi)^n}{n!} = c_n^n |B^n|^2,  \]
where $c_n > 1$ and $\lim_{n\rightarrow \infty} c_n = 1$.

The reason we do not recover the sharp form of Theorem \ref{thm:ineqVolProdMonotone} is the same reason, Urysohn's inequality (\ref{Uryinequ}) is not a special case of its functional analogue (\ref{functurysohn}), namely, that
extremizers in the functional inequalities are Gaussians in both cases while the relevant geometric quantities are recovered for indicators of convex bodies.
However, let us emphasize that the weakest inequality of Theorem \ref{thm:ineqVolProdFunc}, obtained for $\mu = \sigma$, yields a new functional analogue of Urysohn's inequality from which (\ref{Uryinequ}) can be deduced in an asymptotically optimal way, in contrast to (\ref{functurysohn}). We will make this even more precise with our final result which shows that all inequalities of Theorem \ref{thm:ineqVolProdFunc} are strictly stronger than (\ref{functurysohn}). The proof uses ideas from \textbf{\cite{Rotem2012}} and relies on a basic inequality from information theory, sometimes attributed to Shannon (see \textbf{\cite[\textnormal{Theorem~B.1}]{McEliece2002}}), which states that if $g,h: \mathbb{R}^n \rightarrow \mathbb{R}$ are non-negative measurable functions such that $g > 0$ and $\int_{\mathbb{R}^n}g\,dx = 1$, then
\begin{equation} \label{shannoninequ}
\int_{\mathbb{R}^n} g \log \frac{1}{h}\,dx \geq \int_{\mathbb{R}^n} g\log \frac{1}{g}\,dx - \log \left ( \int_{\mathbb{R}^n} h\,dx \right )
\end{equation}
with equality if and only if $h = \alpha g$ for some $\alpha \geq 0$ almost everywhere.

\begin{theorem} Let $\mu$ be an $\mathrm{SO}(n-1)$ invariant probability measure on $\mathbb{S}^{n-1}$ with center of mass at the origin. If $f\in \mathrm{LC_c}(\mathbb{R}^n)$, then
\begin{align} \label{lastthminequ}
\frac{2}{n} \int_{\mathbb{R}^n} h(f,x)\,d\gamma_n(x) \geq 1 - \frac{2}{n} \log\!\left (\! \frac{1}{(2\pi)^{n/2}} \int_{\mathbb{R}^n} (\Psi_\mu f)^\circ(x)\,dx \!\right )
 \end{align}
 with equality if and only if $\Psi_\mu f$ is a multiple of the standard Gaussian.
\end{theorem}

\noindent {\it Proof.} First note that, by Lemma \ref{meanwidthasplundendo},
\[\frac{2}{n}\! \int_{\mathbb{R}^n}\! h(f,x)\,d\gamma_n(x) = \frac{2}{n}\! \int_{\mathbb{R}^n}\! h(\Psi_{\mu}f,x)\,d\gamma_n(x)
= \frac{2}{n}\!\int_{\mathbb{R}^n}\! \log\left(\! \frac{1}{e^{-h(\Psi_\mu f,x)}}\! \right)\! \psi_n(x)\,dx.  \]
Choosing $g = \psi_n$ and $h = e^{-h(\Psi_\mu f,\cdot)}$ in (\ref{shannoninequ}), we thus obtain
\begin{equation} \label{specialshannon}
\frac{2}{n}\! \int_{\mathbb{R}^n}\! h(f,x)\,d\gamma_n(x) \geq \frac{2}{n}\!\int_{\mathbb{R}^n}\!\psi_n(x) \log\!\left(\!\frac{1}{\psi_n(x)}\!\right)\!dx
-\frac{2}{n}\! \log \left( \int_{\mathbb{R}^n}\!e^{-h(\Psi_\mu f,x)}\,dx\! \right)\!.
\end{equation}

\pagebreak

The first integral on the right hand side is the entropy of the standard normal distribution which is well known to be $\frac{n}{2}\left (1+\log(2\pi) \right )$.
Consequently, we obtain
\[\frac{2}{n}\! \int_{\mathbb{R}^n}\! h(f,x)\,d\gamma_n(x) \geq 1 + \log(2\pi)- \frac{2}{n} \log\!\left (\int_{\mathbb{R}^n} (\Psi_\mu f)^\circ(x)\,dx \!\right )\]
which is clearly equivalent to (\ref{lastthminequ}). Equality holds in (\ref{lastthminequ}) if and only if we have equality in (\ref{specialshannon}), that is, by the equality
conditions of (\ref{shannoninequ}), if and only if $ e^{-h(\Psi_\mu f,\cdot)} = \alpha \psi_n,$ for some $\alpha > 0$ or, equivalently, if and only if
\begin{align*}
\mathcal{L}(-\log \Psi_\mu f)(x) = h(\Psi_\mu f,x) = \frac{\|x\|^2}{2} + \beta,
\end{align*}
for some $\beta \in \mathbb{R}$ and every $x \in \mathbb{R}^n$. This shows, by Proposition \ref{propLegendre} (i) and (iv), that equality holds in (\ref{lastthminequ}) if and only if
$\Psi_\mu f = e^{-\beta}(2\pi)^{n/2}\psi_n$. \hfill $\blacksquare$

\vspace{1cm}

\centerline{\large{\bf{ \setcounter{abschnitt}{5}
\arabic{abschnitt}. Appendix}}}

\reseteqn \alpheqn \setcounter{theorem}{0}

\vspace{0.6cm}

The purpose of this appendix is to complete the proof of Proposition \ref{propforthm3} by showing that for $\varphi \in \mathrm{Cvx}_{(o)}(\mathbb{R}^n)$, the function
$\varphi \circledast \mu$ is convex. The proof is based on arguments used in \textbf{\cite{Kiderlen2006}} and \textbf{\cite{Schneider1974c}}, where variants of this fact were
shown under additional assumptions on the function $\varphi$. We begin with an auxiliary result.

\begin{lem} \label{lem:gztausigmaConvex}
Let $\varphi : \mathbb{R}^n \rightarrow \mathbb{R}$ be convex and $H \subseteq \mathbb{R}^n$ a $2$-dimensional linear subspace.
For every $z \in \mathbb{R}^n$, $a, b \in \mathbb{R}$, the function $g_{z, a, b}: H \rightarrow \mathbb{R}$, defined by
\begin{align} \label{defgzab}
g_{z, a, b}(x) = \varphi(a x + b \vartheta_H x + \|x\| z) + \varphi(a x + b \vartheta_H x - \|x\| z),
\end{align}
where $\vartheta_H \in \mathrm{SO}(n)$ acts as rotation by the angle $\frac{\pi}{2}$ on $H$ and keeps $H^\perp$ fixed, is convex.
\end{lem}

\noindent {\it Proof.} Since $\varphi$ and, thus, $g_{z,a,b}$ are continuous, it is sufficient to show that
 \begin{align} \label{whatwewant17}
  g_{z,a,b}\!\left(\frac{x + y}{2}\right) \leq \frac{1}{2}\,g_{z,a,b}(x) + \frac{1}{2}\,g_{z,a,b}(y)
 \end{align}
for all distinct $x, y \in H$. As, by definition, $g_{z,a,b}\!\left(\frac{x + y}{2}\right)$ equals
\begin{align} \label{yoda}
\varphi\!\left(\! a\,\frac{x+y}{2} + b\vartheta_H \frac{x+y}{2}+ \frac{\|x+y\|}{2}z\!\right)\!+ \varphi\!\left(\! a\,\frac{x+y}{2} + b\vartheta_H \frac{x+y}{2}- \frac{\|x+y\|}{2}z\!\right)\!\!,
\end{align}
we may only consider the first term for the following computation and then flip the sign of $z$.
In order to see (\ref{whatwewant17}), first note that for every $\alpha \in [0,1]$,
\begin{align*}
 &a\, \frac{x+y}{2} + b \vartheta_H \frac{x+y}{2}+ \frac{\|x+y\|}{2}z \\
 &= \left(\!a\,\frac{x}{2} + b\vartheta_H \frac{x}{2}+ \alpha \frac{\|x+y\|}{2}z \!\right)
 + \left(\!a\, \frac{y}{2} + b\vartheta_H \frac{y}{2}+ (1-\alpha)\frac{\|x+y\|}{2}z\! \right)\!.
\end{align*}
Again, we may consider only the first term and skip the computation for the second (just replace $x$ by $y$).
Choosing $\alpha = \frac{\|x\|}{\|x\|+\|y\|}$, we have $1-\alpha = \frac{\|y\|}{\|x\| + \|y\|}$ and
\begin{align*}
 &a\,\frac{x}{2} + b\vartheta_H \frac{x}{2}+ \alpha \frac{\|x+y\|}{2}z\\
 &= \lambda_1\! \left(a x + b\vartheta_H x + \|x\|z\right)  +\lambda_2\! \left(a x + b \vartheta_H x - \|x\|z \right)\!,
\end{align*}
where
\begin{align*}
\lambda_1 = \frac{1}{4}\! \left(\!1 + \frac{\|x+y\|}{\|x\| + \|y\|} \right) \qquad \text{ and} \qquad \lambda_2 = \frac{1}{4}\! \left(\!1 - \frac{\|x+y\|}{\|x\| + \|y\|} \right)\!.
\end{align*}
Hence, we conclude that
\begin{align*}
 &a\,\frac{x+y}{2} + b \vartheta_H \frac{x+y}{2}+ \frac{\|x+y\|}{2}z\\
 &= \lambda_1 (a x + b \vartheta_H x + \|x\|z) + \lambda_2 (a x + b \vartheta_H x - \|x\|z) \\
 &\quad + \lambda_1 (a y + b\vartheta_H y + \|y\|z) + \lambda_2 (a y + b \vartheta_H y - \|y\|z).
\end{align*}
Noting that $\lambda_1, \lambda_2 \in \left [0,\frac{1}{2}\right]$ and that $2\lambda_1 + 2 \lambda_2 = 1$, the convexity of $\varphi$ implies that
\begin{align*}
 &\varphi\!\left(\!a\,\frac{x+y}{2} + b\vartheta_H \frac{x+y}{2}+ \frac{\|x+y\|}{2}z\!\right)\\
&\leq \lambda_1 \varphi(a x + b \vartheta_H x + \|x\|z) + \lambda_2 \varphi(a x + b \vartheta_H x - \|x\|z) \\
 &\quad + \lambda_1 \varphi(a y + b\vartheta_H y + \|y\|z) + \lambda_2 \varphi(a y + b \vartheta_H y - \|y\|z).
\end{align*}
The analogue computation for the second term of (\ref{yoda}), finally yields the desired inequality (\ref{whatwewant17}). \hfill $\blacksquare$

\vspace{0.3cm}

In order to show that for  $\varphi \in \mathrm{Cvx}_{(o)}(\mathbb{R}^n)$, the function $\varphi \circledast \mu$ is convex, we first assume that
$\varphi$ is convex and finite on $\mathbb{R}^n$. Moreover, we may restrict ourselves to convex combinations along lines that lie completely in
$\mathbb{R}^n \backslash \{0\}$, the general case then follows by continuity.

Since a zonal function on $\mathbb{S}^{n-1}$ depends only on the value of $u\cdot \bar{e}$, there is a natural one-to-one correspondence between zonal functions and measures on $\mathbb{S}^{n-1}$ and functions and measures on $[-1,1]$ (see, e.g., \textbf{\cite{Schu06a}}). In particular, there exists a unique non-negative measure $\widehat{\mu}$ on $[-1,1]$ such that for every
$f \in C(\mathbb{S}^{n-1})$, we have
\[\int_{\mathbb{S}^{n-1}}\!\!\! f(v)\,d\mu(v) = \int_{-1}^1 \int_{\mathbb{S}^{n-1} \cap \bar e^\perp}\!\!\! f\!\left (\!\alpha \bar{e} + \sqrt{1-\alpha^2} w\!\right) d\sigma_{\bar{e}^\perp}\!(w)\,
(1-\alpha^2)^{\frac{n-3}{2}}d\widehat{\mu}(\alpha),  \]
where $\sigma_{\bar{e}^{\perp}}$ is the invariant probability measure on $\mathbb{S}^{n-1} \cap \bar e^\perp$. Applying this to definition (\ref{defcircledast}), we obtain for $x \in \mathbb{R}^n\backslash\{0\}$,
\begin{align*}
(\varphi \circledast \mu)(x) &= \int_{-1}^1 \int_{\mathbb{S}^{n-1} \cap \bar e^\perp}\!\!\! \varphi\!\left(\!\alpha x +
\sqrt{1-\alpha^2}\|x\| \vartheta_x  w\!\right) d\sigma_{\bar{e}^\perp}\!(w)\,(1-\alpha^2)^{\frac{n-3}{2}}d\widehat\mu(\alpha)\\
&= \int_{-1}^1 \int_{\mathbb{S}^{n-1} \cap x^\perp}\!\!\! \varphi\!\left (\!\alpha x + \sqrt{1-\alpha^2}\|x\| v\!\right) d\sigma_{x^{\perp}}\!(v)\,(1-\alpha^2)^{\frac{n-3}{2}}d\widehat\mu(\alpha).
\end{align*}
Since $\widehat{\mu}$ is non-negative, we are done if we can prove the convexity of the function
\begin{align*}
\varphi_\alpha (x) = \int_{\mathbb{S}^{n-1} \cap x^\perp}\!\!\! \varphi\!\left (\!\alpha x + \sqrt{1-\alpha^2}\|x\| v\!\right) d\sigma_{x^{\perp}}\!(v),
\end{align*}
for all $\alpha \in [-1, 1]$. To this end, let $x \in \mathbb{R}^n \backslash\{0\}$ and let $H \subseteq \mathbb{R}^n$ be an arbitrary $2$-dimensional linear subspace containing $x$. Then $x^\perp = H^\perp \oplus \mathrm{span}\{\vartheta_H x\}$.

First, consider the case $n=3$, where $H = w^\perp$ for some non-zero $w \in \mathbb{R}^3$.
Using cylindrical coordinates $v = \beta \vartheta_H \frac{x}{\|x\|} \pm \sqrt{1-\beta^2}w$ on $\mathbb{S}^2 \cap x^\perp$, we obtain
\begin{align*}
 \varphi_\alpha(x) &= \frac{1}{\pi}\!\int_{-1}^1\! \varphi\!\left(\!\alpha x + \sqrt{1-\alpha^2}\beta \vartheta_H x + \sqrt{(1-\alpha^2)(1-\beta^2)}\|x\|w\! \right)\! \frac{d \beta}{\sqrt{1-\beta^2}} \\
                    & \phantom{=} + \frac{1}{\pi}\!\int_{-1}^1\! \varphi\!\left(\!\alpha x + \sqrt{1-\alpha^2}\beta \vartheta_H x - \sqrt{(1-\alpha^2)(1-\beta^2)}\|x\|w\! \right)\! \frac{d \beta}{\sqrt{1-\beta^2}}\\
                   &= \frac{1}{\pi}\!\int_{-1}^1\! g_{z,a,b}(x) \frac{d \beta}{\sqrt{1-\beta^2}},
\end{align*}
where $z = \sqrt{(1-\alpha^2)(1-\beta^2)}w$, $a = \alpha$ and $b = \sqrt{1-\alpha^2}\beta$. By Lemma \ref{lem:gztausigmaConvex}, $g_{z,a,b}$ is convex and, hence, $\varphi_\alpha$ is convex, as well.

For $n \geq 4$, we again use cylindrical coordinates on $\mathbb{S}^{n-1} \cap x^\perp$ in the direction of $\vartheta_H \frac{x}{\|x\|}$ to obtain
\begin{align*}
 \varphi_\alpha(x) = c_n\!\! \int_{-1}^1\! \int_{\mathbb{S}^{n-1}\cap H^\perp}\hspace{-1.2cm} \varphi ( \alpha x\! +\! \sqrt{1\!-\!\alpha^2}\beta \vartheta_H x\! +\! \sqrt{(1\!-\!\alpha^2)(1\!-\!\beta^2)}\|x\|w) d \sigma_{\!H^\perp\!}(v) \frac{d \beta}{(1\!-\!\beta^2)^\frac{4-n}{2}}\! ,
\end{align*}
where $c_n = \frac{2\Gamma(\frac{n-1}{2})}{\sqrt{\pi}\Gamma(\frac{n}{2}-1)}$. Taking $\frac{1}{2}$ of this integral twice and replacing $v$ by $-v$ in one copy, we see that again
\begin{align*}
 \varphi_\alpha(x) = \frac{c_n}{2}\! \int_{-1}^1\! \int_{\mathbb{S}^{n-1}\cap H^\perp} \hspace{-0.3cm} g_{z,a,b}(x)\, d\sigma_{\!H^\perp\!}(v) \frac{d \beta}{(1\!-\!\beta^2)^\frac{4-n}{2}}\!,
\end{align*}
where $z = \sqrt{(1-\alpha^2)(1-\beta^2)}w$, $a = \alpha$ and $b = \sqrt{1-\alpha^2}\beta$ as before. By Lemma~\ref{lem:gztausigmaConvex}, $g_{z,a,b}$ is convex and, thus, so is $\varphi_\alpha$.

For general $\varphi \in \mathrm{Cvx}_{(o)}(\mathbb{R}^n)$, we use that the Moreau envelope $e_t\varphi$ of $\varphi$ is convex and finite and converges monotonously to $\varphi$,
by Lemma \ref{propMoreau}. By what we have shown above, each of the functions $e_t\varphi \circledast \mu$, $t > 0$, is convex. Hence, by monotone convergence, we
conclude that
 \begin{align*}
  \varphi \circledast \mu = \lim_{t \searrow 0} e_t\varphi \circledast \mu = \sup_{t > 0} e_t\varphi \circledast \mu
 \end{align*}
is convex as well.

\vspace{0.6cm}

\noindent {{\bf Acknowledgments} The authors were supported by the Austrian Science Fund (FWF), Project number: P31448-N35.

\begin{small}

\[ \begin{array}{ll} \mbox{Georg C. Hofst\"atter} & \mbox{Franz E. Schuster} \\
\mbox{Vienna University of Technology \phantom{wwwwWW}} & \mbox{Vienna University of Technology} \\
\mbox{Inst.\ for Discrete Math.\ and Geometry} & \mbox{Inst.\ for Discrete Math.\ and Geometry} \\
\mbox{Wiedner Hauptstrasse 8-10/1047} & \mbox{Wiedner Hauptstrasse 8-10/1047} \\
\mbox{1040 Vienna, Austria} & \mbox{1040 Vienna, Austria} \\
\mbox{georg.hofstaetter@tuwien.ac.at} & \mbox{franz.schuster@tuwien.ac.at}
\end{array}\]

\end{small}

\end{document}